\RequirePackage[l2tabu,orthodox]{nag}	%LaTeXe orthodoxe ?

\documentclass[a4paper,12pt,english]{article}

	\usepackage{lmodern}				%Latin modern, police vectorielle

	\usepackage[T1]{fontenc}			%Occidental

	\usepackage[utf8x]{inputenc}		%UTF8

	\usepackage{babel}					%Typographie

	\usepackage{amsmath}		%Math de l'AMS
	\usepackage{amssymb}		%Autres symboles mathématiques
	\usepackage{amsthm}			%Théorèmes correctement organisés

	\usepackage{graphicx}	%Insérer des images
	\usepackage{color}		%Couleurs (nécessaires pour les images avec XFig)

	\usepackage{version}		%Commenter (notamment un environnement)

	\usepackage{geometry}		%Marges

\input{commath}

	\usepackage{hyperref}				%Hyperliens
	
	\usepackage{cleveref}				%Références intelligentes

	\usepackage{pstricks}
	\usepackage{ifmtarg}
	\usepackage{tikz,tkz-tab}
	\usepackage{authblk}
	
\makeatletter
\newcommand*{\rva}[3][]{%
	\@ifmtarg{#1}{%		#1 est-il vide ?
		#2_{#3}%			oui : X_n
	}{%
		#2_{#3,#1}%		non : X_{n,k}
	}%
}
\makeatother

\newcommand*{\Xinf}{X^<}
\newcommand*{\Sinf}{S^<}
\newcommand*{\Yinf}{Y^<}
\newcommand*{\Tinf}{T^<}
\newcommand*{\Xtilde}{\tilde{X}}
\newcommand*{\Ytilde}{\tilde{Y}}

\newcommand*{\ssdd}[1]{\sigma_{#1}}

\renewcommand*{\th}{th}

\theoremstyle{plain}
\newtheorem{thm}{Theorem}

\newtheorem{lem}[thm]{Lemma}
\newtheorem{prop}[thm]{Proposition}

\theoremstyle{remark}
\newtheorem{rem}[thm]{Remark}

\newcommand{\pscal}[1]{\langle #1 \rangle}

\DeclareMathOperator{\Hess}{Hess}

\hypersetup{%
colorlinks=true,%
linkcolor = black,%
urlcolor = blue,%
citecolor = blue}

\makeatletter
\def\input@path{{./images/}}
\makeatother

\begin{document}

\title{Deviation results for sparse tables in hashing with linear probing}

\author[1]{Thierry Klein}
\author[2]{Agn\`es Lagnoux}
\author[3]{Pierre Petit}

\affil[1]{Institut de Math\'ematiques de Toulouse; UMR5219. Universit\'e de Toulouse; ENAC - Ecole Nationale de l'Aviation Civile, Universit\'e de Toulouse, France. E-mail: thierry.klein@math.univ-toulouse.fr}
\affil[2]{Institut de Math\'ematiques de Toulouse; UMR5219. Universit\'e de Toulouse; CNRS. UT2J, F-31058 Toulouse, France. E-mail: lagnoux@math.univ-tlse2.fr (corresponding author)}
\affil[3]{Institut de Math\'ematiques de Toulouse; UMR5219. Universit\'e de Toulouse; CNRS. UT3, F-31062 Toulouse, France. E-mail: pierre.petit@math.univ-toulouse.fr}

\date{12 November 2021}

\maketitle

\begin{abstract}
We consider the model of hashing with linear probing and we establish the moderate and large deviations for the total displacement in sparse tables. In this context, Weibull-like-tailed random variables  appear. Deviations for sums of such heavy-tailed random variables  are studied in \cite{Nagaev69-1,Nagaev69-2}.  Here we adapt the proofs therein to deal with conditioned sums of such variables and solve the open question in  \cite{TFC12}. By the way, we establish the deviations of the total displacement in full tables, which can be derived from the deviations of empirical processes of i.i.d.\ random variables established in \cite{Wu94}.
\end{abstract}



\section{Introduction}

Hashing with linear probing is a classical model in theoretical computer science that appeared in the 50's. It has been studied from a mathematical point of view firstly by Knuth in \cite{Knuth63}. 
Here is a simple description given in \cite{FPV98}.
\begin{quote}
A table of length $m$, $T[1.. m]$ is set up, as well as a hash function $h$ that maps
keys from some domain to the interval $[1.. m]$ of table addresses. A collection
of $n$ elements with $n\leq m$ are entered sequentially into the table according to the
following rule: Each element $x$ is placed at the first unoccupied location starting from
$h(x)$ in cyclic order, namely the first of $h(x)$, $h(x)+1$, \ldots, $m$, $1$, $2$, \ldots, $h(x)-1$.
\end{quote}
For more details on the model, we refer to \cite{FPV98, Janson01a, Marckert01-1, Chassaing02, ChF03, Janson05}. The length of the move of each element $x$ is called the \emph{displacement} of $x$ and the sum of all displacements, denoted by $d_{m,n}$, is called the \emph{total displacement}. 
In its seminal papers \cite{Knuth63,Knuth74}, Knuth assumes that all the sequence of hash addresses $h(x)$ are independent and uniformly distributed on $\intervallentff{1}{m}$, computes exact expressions of $\Espe[d_{m,n}]$ and $\Var(d_{m,n})$, and provides their asymptotic behaviors. The limit distribution of $d_{m,n}$ remains unknown until 1998: in \cite{FPV98}, Flajolet, Poblete, and Viola give the limit distribution of $d_{m,n}$ for full tables ($n=m$) and for sparse tables ($n/m=\mu \in \intervalleoo{0}{1}$) using combinatorial arguments. In \cite{Marckert01-2}, Chassaing and Marckert recover the previous results in the full case via a probabilistic approach. They prove that $d_{m,n}$ is the area under the {\L}ukasiewicz random walk (also called Breadth First Search random walk) associated to a Galton-Watson tree with Poisson progeny. Consequently, the limit distribution of the total displacement $d_{m,m}$ is that of the area under the Brownian excursion, which involves the Airy distribution.

In \cite{Janson01a}, reformulating the problem in terms of conditioned sums of random variables, Janson establishes the limit distribution of $d_{m,n}$ in all cases with probabilistic tools. 
In \cite{Janson01}, Janson extends the central limit theorem in the sparse case to a general model of conditioned sums of random variables. The corresponding Berry-Esseen bounds are proved by Klein, Lagnoux, and Petit in \cite{TAP19}. 
Concerning the deviations of such conditioned models, Gamboa, Klein, and Prieur give an answer in the case of light-tailed random variables (see \cite{TFC12}). Unfortunately, their results cannot be applied to the model of hashing with linear probing since this model involves heavy-tailed random variables.  

In this paper, we establish the moderate and large deviations for the total  displacement $d_{m,n}$ in sparse tables. Deviations for heavy-tailed random variables are studied by several authors (e.g., \cite{Petrov54, Linnik61, Nagaev69-1, Nagaev69-2, Nagaev_1979_AnnProb}) and a good survey can be found in Mikosch \cite{mikosch1998large}. 
In the context of hashing with linear probing, Weibull-like-tailed random variables appear. Deviations for sums of such variables  are studied  by Nagaev in \cite{Nagaev69-1,Nagaev69-2}.  Here we adapt its proofs  to deal with conditioned sums of such variables. We also need to establish the deviations of $d_{m,n}$ for full tables, which can be derived from the deviations of empirical processes of i.i.d.\ random variables established by Wu in \cite{Wu94}.

The paper is organized as follows. In Section \ref{sec:main}, we state the main results for full and sparse tables. The proofs for full tables are given in Section  \ref{sec:proof_full} and those for sparse tables can be found in Section \ref{sec:proof_sparse}. In Section \ref{sec:inter}, we expose Janson's reformulation of the model and provide several useful estimates required in Section \ref{sec:proof_sparse}.

\section{Setting and main results}\label{sec:main}

\subsection{Model}

An equivalent formulation of the problem of hashing can be made in terms of the  discrete version of the classical parking problem described, for instance, by Knuth \cite{Knuth98}:

\begin{quote}
A certain one-way street has m parking spaces in a row numbered $1$ to $m$. A
man and his dozing wife drive by, and suddenly, she wakes up and orders him to park
immediately. He dutifully parks at the first available space [...]. 
\end{quote}

More precisely, the model describes the following experiment. Let $n \leqslant m$. $n$ cars enter sequentially into a circular parking uniformly at random. 
The parking spaces are numbered clockwise. A car that intends to park at an occupied space moves to the next empty space, always moving clockwise. 
The length of the move is called the displacement of the car and we are interested in the sum of all displacements which is a random variable denoted by $d_{m,n}$. When all cars are parked, there are $N = m-n$ empty spaces. 
These divide the occupied spaces into blocks of consecutive spaces. We consider that the empty space following a block belongs to this block.

For example, assume that $n=8$, $m=10$, and $(6,9,1,9,9,6,2,5)$ are the addresses where the cars land. This sequence of (hash) addresses is called a \emph{hash sequence} of length $m$ and size $n$. Let $d_i$ be the displacement of car $i$. Then $d_1=d_2=d_3=0$. The car number $4$ should park on the 9\th space which is occupied by the 2nd car; thus it moves one space ahead and parks on the 10\th space so that $d_4=1$. The car number $5$ should park on the 9\th space. Since the 9\th, the 10\th, and the 1st spaces are occupied, $d_5=3$. And so on: $d_6=1,\ d_7=1,\ d_8=0$. Here, the total displacement is equal to $d_{10,8}=1+3+1+1=6$. In our example, there are two blocks: the first one containing spaces $9$, $10$, $1$, $2$, $3$ (occupied), and space $4$ (empty), and the second one containing spaces $5$, $6$, $7$ (occupied), and space $8$ (empty).

In this paper, we are interested in the deviations of the total displacement $d_{m,n}$ in sparse tables. To do so, we need the large deviation behavior of $d_{m,n}$ in full tables. By the way, we also established the moderate deviation behavior of $d_{m,n}$ in full tables.

\subsection{Deviations in full tables}

In this section, we first recall some already existing results for the total displacement $d_{m,m}$ in full tables ($n=m$). As mentioned in the introduction, Knuth in \cite{Knuth74} and Flajolet et al.\ in \cite[Theorem 2]{FPV98}) derive the asymptotic behavior of the expectation and the variance of $d_{m,m}$: 
\begin{align}\label{eq:behav_moments_dense}
\Espe[d_{m,m}] \underset{m\to \infty}{\sim} \frac{\sqrt{2\pi}}{4} m^{3/2}
\quad \text{and} \quad
\Var(d_{m,m}) \underset{m\to \infty}{\sim} \frac{10-3\pi}{24} m^{3}.
\end{align}

The following result was first established in \cite[Theorem 3]{FPV98}. 

\begin{thm}[Standard deviations] \label{th:sd_full} 
For full tables, the distribution of the total displacement $d_{m,m}/m^{3/2}$ is asymptotically distributed as the area $A$ under the standard Brownian excursion, in the sense that, for all $\delta\geqslant 0$,
\[
\Prob(d_{m,m} \geqslant m^{3/2} \delta) \xrightarrow[m\to \infty]{} \Prob(A\geqslant \delta).
\]
\end{thm}

In this paper, we establish the probabilities of deviation for the total displacement $d_{m,m}$.

\begin{thm}[Moderate deviations]\label{th:md_full}
For all $\alpha \in \intervalleoo{3/2}{2}$ and for all $\delta\geqslant 0$, 
\[
\frac{1}{m^{2\alpha-3}} \log \Prob(d_{m,m} \geqslant m^\alpha \delta) \xrightarrow[m\to \infty]{} -6 \delta^2.
\]
\end{thm}

\begin{thm}[Large deviations] \label{th:ld_full}
For all $\delta\geqslant 0$, 
\[
- \frac{1}{m} \log \Prob(d_{m,m} \geqslant m^2 \delta)
 \xrightarrow[m\to \infty]{} J(\delta)
 \defeq \begin{cases}
(\frac{1}{2}-\delta) \cdot{} \lambda(\delta) +\log(1-(\frac{1}{2}+\delta)\cdot{} \lambda(\delta)) & \text{if $\delta<1/2$}\\
\infty & \text{if $\delta\geqslant 1/2$,}
\end{cases}
\]
where $\lambda(\delta)$ is the smallest solution of the equation in $\lambda$ 
\begin{equation}\label{eq:lambda}
\left(\lambda \cdot \left(\delta+\frac 12\right)-1\right)(1-e^\lambda) = \lambda .
\end{equation}
\end{thm}

Observe that lower deviations are trivial: for all $\alpha \in \intervalleof{3/2}{2}$ and all $m$ large enough,
$\Prob(d_{m,m}-\Espe[d_{m,m}] \leqslant - m^{\alpha} \delta) = 0$
because of the positiveness of $d_{m,m}$ and using \eqref{eq:behav_moments_dense}. When dealing with the very large deviations, the same trivial behavior occurs both for upper and lower deviations: for all $\alpha > 2$ and all $m$ large enough,
$
\Prob(d_{m,m}-\Espe[d_{m,m}] \geqslant m^{\alpha} \delta) = 0$ and $
\Prob(d_{m,m}-\Espe[d_{m,m}] \leqslant - m^{\alpha} \delta) = 0$,
since $d_{m,m} \leqslant m(m-1)/2$.

\begin{rem}
If $\delta=0$, $\lambda(0)=0$ is the unique solution of Equation \eqref{eq:lambda}. If $\delta\in \intervalleoo{0}{1/2}$, Equation \eqref{eq:lambda} has two solutions: $\lambda(\delta)<0$ and $0$. 
Moreover, $J(\delta)\sim 6\delta^2$ as $\delta\to 0$ (we recover the rate function of the moderate deviations in Theorem \ref{th:md_full}) and $J(\delta)\to +\infty$ as $\delta\to 1/2$.
\end{rem}

\begin{rem}\label{rem:dense_case}
The conclusions of Theorems \ref{th:sd_full}, \ref{th:md_full}, and \ref{th:ld_full} are still valid replacing $\delta$ by $\delta + o(1)$, since the limiting functions are continuous. In particular, one may replace $d_{m,m}$ by $d_{m,n}$ as soon as $m-n \ll m^{\alpha-1}$ (for instance, in the almost full case where $n=m-1$). Indeed, naturally coupling $d_{m,n}$ and $d_{m,m}$ by adding $m-n$ balls, one has
\begin{align*}
\abs{d_{m,n}-d_{m,m}}
 & \leqslant (m-1) + (m-2) + \dots + n
% & = \frac{(m+n-1)(m-n)}{2}\\
 \sim m (m-n)
 \ll m^{\alpha} ,
\end{align*}
whence
\[
\Prob(d_{m,n} \geqslant m^\alpha \delta)
 \leqslant \Prob\Bigl( d_{m,m} \geqslant m^\alpha \Bigl( \delta - \frac{m(m-n)}{m^\alpha} \Bigr) \Bigr)
 = \Prob(d_{m,m} \geqslant m^\alpha(\delta + o(1)))
\]
and, similarly,
\[
\Prob(d_{m,n} \geqslant m^\alpha \delta)
 \geqslant \Prob\Bigl( d_{m,m} \geqslant m^\alpha \Bigl( \delta + \frac{m(m-n)}{m^\alpha} \Bigr) \Bigr)
 = \Prob(d_{m,m} \geqslant m^\alpha(\delta + o(1))) .
\]
%\[
%\Prob(d_{m,n} \geqslant m^\alpha \delta) = \Prob(d_{m,m} \geqslant m^\alpha(\delta + o(1))) .
%\]
%So our results complement the standard deviations established in \cite{Janson01a}, in the dense case. 
Moreover, using a new probabilistic approach developed in \cite{Marckert01-2}, Theorem \ref{th:sd_full} was extended in \cite[Theorems 1.1 and 2.2]{Janson01a} to the case $(m-n)/\sqrt{m} \to a \in \intervallefo{0}{\infty}$: for all $\delta\geqslant 0$,
\[
\Prob(d_{m,n} \geqslant m^{3/2} \delta) \xrightarrow[m\to \infty]{} \Prob(W_a\geqslant \delta),
\]
with 
\[
W_a=\int_0^1 \max_{s\leqslant t} (b(t)-b(s)-a(t-s))dt
\]
where $b$ is a Brownian bridge $b$ on $\intervalleff{0}{1}$ periodically extended to $\R$.
\end{rem}

\subsection{Deviations in sparse tables}

In this section, we consider asymptotics in $(m,n)$ with $m \to \infty$ and $n/m \to \mu \in \intervalleoo{0}{1}$ (sparse case).  This definition of the sparse case is a slight extension of that of \cite{FPV98} ($n/m = \mu \in \intervalleoo{0}{1}$). Set $N = m-n$.
%; but it is worth mentioning that the sparse case considered by Janson in \cite{Janson01a} has a wider range: $N \gg \sqrt{m}$ and $n \gg \sqrt{m}$.
%Notice that $m = \Theta(n)$ and $N = m-n = \Theta(m)$.
%Here, it is natural to make quantities depend on the number $N$ of empty urns. Let $(n_N)_{n \geqslant 1}$ be a sequence of positive integers and let $m_N= n_N+N$. Assume that $n_N/m_N \to \mu \in \intervalleoo{0}{1}$ as $N\to \infty$ (sparse case).
By \cite[Theorem 5]{FPV98},
\begin{align}\label{eq:behav_moments_full}
\Espe[d_{m,n}] \underset{m\to \infty}{\sim} \frac{\mu^2}{2(1-\mu)^2} N \quad \text{and} \quad 
\Var(d_{m,n}) \underset{m\to \infty}{\sim} \sigma^2(\mu)N,
\end{align}
where (cf. \cite[Theorem 5]{FPV98})
\begin{align}\label{def:sigma_mu}
\sigma^2(\mu) \defeq \frac{6\mu^2 - 6\mu^3 + 4\mu^4 - \mu^5}{12(1-\mu)^5} .
\end{align}

The following result was first proved in \cite[Theorem 6]{FPV98} while another probabilistic proof was given in \cite{Janson01a}.

\begin{thm}[Standard deviations] \label{th:sd_sparse}
The distribution of the total displacement $d_{m,n}$ is asymptotically Gaussian distributed, in the sense that, for all $y$,
\[
\Prob\left(d_{m,n}-\Espe[d_{m,n}] \leqslant N^{1/2} y\right) \xrightarrow[m\to \infty]{} \Prob(Z \leqslant y)
\]
where $Z \sim \mathcal{N}(0,\sigma^2(\mu))$. 
\end{thm}

In this paper, we establish the probabilities of deviation of the total displacement $d_{m,n}$ in the sparse case.

\begin{thm}[Lower moderate deviations] \label{th:lmd_sparse}
For all $\alpha\in \intervalleoo{1/2}{1}$ and for all $y \geqslant 0$,
\begin{align} \label{eq:lmd_sparse}
\frac{1}{N^{2\alpha-1}} \log \Prob(d_{m,n} - \Espe[d_{m,n}] \leqslant - N^{\alpha} y)
 & \xrightarrow[m\to \infty]{} - \frac{y^2}{2\sigma^2(\mu)} .
\end{align}
\end{thm}

\begin{thm}[Lower large deviations] \label{th:lld_sparse}
For all $y \geqslant 0$,
\begin{align} \label{eq:lld_sparse}
\frac{1}{N} \log \Prob(d_{m,n} - \Espe[d_{m,n}] \leqslant -N y) \xrightarrow[m\to \infty]{} -\Lambda^*\Bigl( \frac{1}{1-\mu}, \frac{\mu^2}{2(1-\mu)^2}-y \Bigr) ,
\end{align}
where $\Lambda^*$ is the Fenchel-Legendre transform of the function $\Lambda\colon \R^2\to \intervalleof{-\infty}{\infty}$ defined by
\[
\Lambda(s,t) = \log \sum_{l=1}^\infty \frac{e^{(s-\mu)l} (\mu l)^{l-1}}{l!} \Espe[e^{t d_{l,l-1}}] .
\]
%\alert{Expression plus simple ??? Flajolet ne donne que des expressions des moments d'ordre 1 et 2 des almost full tables...}
\end{thm}

For all $\alpha>1$, we have, asymptotically,
$\Prob(d_{m,n} - \Espe[d_{m,n}] \leqslant -N^\alpha y) = 0$, since $d_{m,n} \geqslant 0$ and $\Espe[d_{m,n}]$ is asymptotically linear in $N$.

\begin{thm}[Upper deviations] \label{th:umld_sparse} \leavevmode

(i) For all $\alpha\in \intervalleoo{1/2}{2/3}$ and for all $y \geqslant 0$,
\begin{align} \label{eq:umd_sparse}
\frac{1}{N^{2\alpha-1}}\log \Prob(d_{m,n} - \Espe[d_{m,n}] \geqslant N^{\alpha} y) \xrightarrow[m\to \infty]{} -\frac{y^2}{2\sigma^2(\mu)}.
\end{align}

(ii) For all $y \geqslant 0$,
\begin{align}\label{eq:uid_sparse}
\frac{1}{N^{1/3}} \log \Prob(d_{m,n} - \Espe[d_{m,n}] \geqslant N^{2/3} y) \xrightarrow[m\to \infty]{} - I(y)
\end{align}
with 
\[
I(y)
 \defeq \begin{cases}
\frac{y^2}{2\sigma^2(\mu)} & \text{if $y\leqslant y(\mu)$},\\ 
q(\mu)(1-t(y))^{1/2} y^{1/2}+\frac{t(y)^2y^2}{2\sigma^2(\mu)} & \text{if $y> y(\mu)$}, 
\end{cases}
\]
where 
\[
q(\mu) \defeq \inf_{0<\delta<1/2} \frac{1}{\sqrt{\delta}}(\kappa(\mu) + J(\delta) ),
\]
$\kappa(\mu)\defeq \mu-\log(\mu)-1\in \intervalleoo{0}{\infty}$, $J$ has been defined in Theorem \ref{th:ld_full},
 $y(\mu)\defeq 3\left(q(\mu)\sigma^2(\mu)\right)^{2/3}/2$, and  $t(y)$ is defined for $y>y(\mu)$ as the smallest root of the cubic equation in $t\in \intervalleff{0}{1}$
\begin{align*}
t^3-t^2+\frac{q^2(\mu)\sigma^4(\mu)}{4y^{3}}=0.
\end{align*}

(iii) For all $\alpha \in \intervalleoo{2/3}{2}$ and for all $y \geqslant 0$,
\begin{align}\label{eq:uld_sparse}
\frac{1}{N^{\alpha/2}} \log \Prob(d_{m,n} - \Espe[d_{m,n}] \geqslant N^\alpha y) \xrightarrow[m\to \infty]{} - q(\mu)y^{1/2} .
\end{align}

(iv) For all $y \geqslant 0$,
\begin{align}\label{eq:uld_sparse_particulier}
\frac{1}{N} \log &\ \Prob(d_{m,n} - \Espe[d_{m,n}] \geqslant N^2 y)\nonumber\\ 
&\xrightarrow[m\to \infty]{} \begin{cases}
 - \inf\limits_{\delta > 0} \bigl[ \sqrt{\frac{y}{\delta}}(\kappa(\mu) + J(\delta)) + \Lambda_0^*\bigl(\frac{1}{1-\mu}-\sqrt{\frac{y}{\delta}}\bigr) \bigr]
& \text{if $y<\mu^2/(2(1-\mu)^2)$}\\
-\infty & \text{if $y\geqslant \mu^2/(2(1-\mu)^2)$,}
\end{cases}
\end{align}
where $\Lambda_0^*$ is the Fenchel-Legendre transform of the function $\Lambda_0\colon \R\to \intervalleof{-\infty}{\infty}$ defined by $\Lambda_0(s)\defeq \Lambda(s,0)$ and $\Lambda$ has been defined in Theorem \ref{th:lld_sparse}.
\end{thm}

For all $\alpha > 2$, we have, asymptotically,
$\Prob(d_{m,n} - \Espe[d_{m,n}] \geqslant N^\alpha y) = 0$,
since $d_{m,n} \leqslant n(n-1)/2$ and $N$ is asymptotically linear in $n$.

\begin{rem}
Observe that, for $\alpha=2/3$ and for all $y>0$, 
\[
I(y)=\inf_{t\in \intervalleff{0}{1}} f(t),
\]
where 
\[ 
f(t)=\Bigl(q(\mu)(1-t)^{1/2} y^{1/2}+\frac{t^2y^2}{2\sigma^2(\mu)}\Bigr).
\]

If $y\leqslant y_1(\mu)=3\left(q(\mu)\sigma^2(\mu)\right)^{2/3}/2^{4/3}$, then $f$ is decreasing and its minimum $y^2/(2\sigma^2(\mu))$ is attained at $t=1$. If $y>y_1(\mu)$,
 $f$ has two local minima, at $1$ and at $t(y)$, corresponding to the smallest of the two roots in $\intervalleff{0}{1}$ of $f'(t)=0$ which is equivalent to 
$t^3-t^2+q^2(\mu)\sigma^4(\mu)/(4y^{3})=0$. If $y_1(\mu)< y\leqslant y(\mu)$, the minimum is attained at $1$, and at $t(y)$ otherwise.  
Let $c=q^2(\mu)\sigma^4(\mu)/(4y^{3})$. One can prove that $t(y)=2\Re(z)+1/3$ where $z$ is the only complex cube root  of 
\[
\frac{1}{27}-\frac c2 + i \sqrt{\frac{c}{27}-\frac{c^2}{4}}
\] 
having argument in $\intervalleoo{\pi/3}{2\pi/3}$.
\end{rem}

%\alert{Cas frontière non traité (entre dense et sparse) : déviations en $n^\alpha$ pour $d_{m,n}$ avec $c n^{\alpha-1} \leqslant m-n \ll n$.}

\begin{rem}
One can deduce from the proofs that the following (probably typical) events roughly realize the large deviations in Theorem \ref{th:umld_sparse}:\newline
(i) all the displacements within the blocks are small but their sum has a Gaussian contribution;\newline
(iii) one block has a large size close to $\delta^{-1/2} N^{\alpha/2}y^{1/2}$ and the displacement within this block is close to $N^{\alpha}y$, $\delta$ being chosen by the optimization in $q(\mu)$ (two competing terms);\newline
(ii) one block has a large size and the displacement within this block is large (close to $N_n^\alpha (1-t(y)) y$) and the sum of the other displacements has a Gaussian contribution (two extra competing terms);\newline
(iv) one block has a large size close to $\delta^{-1/2} N y^{1/2}$ and the displacement within this block is close to $N^2y$, which forces the sum of the length of the other blocks to be abnormally small, that is close to $m - \delta^{-1/2} N y^{1/2} \sim N((1-\mu)^{-1} - \delta^{-1/2} y^{1/2})$ (three competing terms).
\end{rem}

%%%%%%%%%%%%%%%%%%%%%%%%%%%%%%%%%%
%%%%%%%%%%%%%%%%%%%%%%%%%%%%%%%%%%
%%%%%%%%%%%%%%%%%%%%%%%%%%%%%%%%%%
\section{Proofs: full tables}\label{sec:proof_full}
%%%%%%%%%%%%%%%%%%%%%%%%%%%%%%%%%%
%%%%%%%%%%%%%%%%%%%%%%%%%%%%%%%%%%
%%%%%%%%%%%%%%%%%%%%%%%%%%%%%%%%%%

Here, we take $n=m$.  All limits are considered as $m\to \infty$ unless stated otherwise.
%\begin{proof}[Proof of Theorem \ref{th:md_full}]
%For $m \geqslant 1$, let $(V_{m,i})_{1 \leqslant i\leqslant m}$ be a sequence of i.i.d.\ random variables uniformly distributed on $\intervallentff{1}{m}$. $V_{m,i}$ corresponds to the hash address of item $i$. Then we define
%\[
%S_m(k)=\sum_{i=1}^m \indic_{V_{m,i}\leqslant k}.
%\]
%
%Now, let $(U_i)_{1\leqslant i\leqslant m}$ be a sequence of i.i.d.\ random variables with uniform distribution $\mathcal{U}$ on $\intervalleff{0}{1}$ and such that, for all $i\in \intervallentff{1}{m}$, $V_{m,i}=\lceil m U_i \rceil$ almost surely. Following \cite{Wu94}, we introduce the empirical measure $L_m$ associated to the random variables $U_i$.
%For $m \geqslant 1$, let $(V_{m,i})_{1 \leqslant i\leqslant m}$ be a sequence of i.i.d.\ random variables uniformly distributed on $\intervallentff{1}{m}$. $V_{m,i}$ corresponds to the hash address of item $i$. Then we define
%\[
%S_m(k)=\sum_{i=1}^m \indic_{V_{m,i}\leqslant k}.
%\]
For all $m\geqslant 1$, let $(U_{m,i})_{1\leqslant i\leqslant m}$ be a sequence of i.i.d.\ random variables with uniform distribution $\mathcal{U}$ on $\intervalleff{0}{1}$ and let, for all $i\in \intervallentff{1}{m}$, $V_{m,i}=\lceil m U_{m,i} \rceil$.
Note that $(V_{m,i})_{1 \leqslant i\leqslant m}$ is a sequence of i.i.d.\ random variables uniformly distributed on $\intervallentff{1}{m}$ and, for all $i=1,\dots,m$, $V_{m,i}$ corresponds to the hash address of item $i$.
For all $k \in \intervallentff{1}{m}$, let us define
\[
S_m(k) \defeq \sum_{i=1}^m \indic_{1 \leqslant V_{m,i} \leqslant k}
 = m L_m(\intervalleff{0}{k/m}),
\]
where $L_m$ is the empirical measure  associated to the random sequence $(U_{m,i})_{1\leqslant i\leqslant m}$. Note that $S_m(m) = m$. As in \cite[Lemma 2.1]{Janson01a}, we extend the definition of $S_m(k)$ for $k \in \Z$ so that the sequence $(S_m(k) - k)_{k \in \Z}$ be $m$-periodic, whence, for all $k \in \Z$,
\[
\min_{l<k} \{S_m(l)-l\} = \min_{l \in \Z} \{S_m(l)-l\} = \min_{1 \leqslant l \leqslant m} \{ S_m(l) - l \} .
\]
Thus, by \cite[Equation (2.1) and Lemma 2.1]{Janson01a}, the total displacement $d_{m,m}$ is given by 
\begin{align} \label{eq:dmm}
d_{m,m}
 & = \biggl( \sum_{k=1}^m S_m(k) - k \biggr) - m \min_{l \in \Z} \{ S_m(l) - l \} .
\end{align}

Since $\min \ensavec{S_m(l) - l}{l \in \Z} \leqslant S_m(m) - m = 0$, for all $z \geqslant 0$,
\begin{align} \label{eq:full_dmm_lb}
\Prob(d_{m, m} \geqslant z)
 & \geqslant \Prob\biggl( \sum_{k=1}^m S_m(k) - k \geqslant z \biggr) .
\end{align}

Let us find an upper bound for the probability in the left-hand side. For all $j \in \intervallentff{1}{m}$, we introduce the sequence $U^j_m \defeq (U_{m,i} - j/m)_{1 \leqslant i \leqslant m}$ where the addition is considered on the torus $\R / \Z$. We also define the associated random variables $V^j_{m, i}$, $S^j_m(k)$, and $d^j_{m, m}$.  The following lemma is straightforward.

\begin{lemma} \label{lem:rotation}
Let $j \in \intervallentff{1}{m}$.

(i) $U^j_m$ has the same distribution as $U^m_m = (U_{m,i})_{1 \leqslant i \leqslant m}$. As a consequence, $((S^j_m(k))_{1 \leqslant k \leqslant m}, d^j_{m, m})$ has the same distribution as $((S_m(k))_{1 \leqslant k \leqslant m}, d_{m, m})$.

% (ii) For all $i \in \intervallentff{1}{m}$, $V^j_{m, i} = V_{m, i} - j$, where the addition is considered on $\Z / (m \Z)$.

(ii) For all $k \in \Z$, $S^j_m(k) = S_m(j+k) - S_m(j)$.

(iii) $d^j_{m, m} = d_{m, m}$.
\end{lemma}

Let $j_0 \in \intervallentff{1}{m}$ be such that
\[
S_m(j_0) - j_0 = \min_{l \in \Z} \{ S_m(l) - l \} .
\]
We claim that
\[
\min_{l \in \Z} \{ S^{j_0}_m(l) - l \} = 0 .
\]
Since $S^{j_0}_m(0) = 0$, it is enough to show that, for all $l \in \Z$, $S^{j_0}_m(l) - l \geqslant 0$. Using Lemma \ref{lem:rotation} (ii),
\[
S^{j_0}_m(l) - l
 = S_m(j_0 + l) - S_m(j_0) - l
 = S_m(j_0 + l) - (j_0 + l) - (S_m(j_0) - j_0)
 \geqslant 0 .
\]

Therefore, for all $z \geqslant 0$, using Lemma \ref{lem:rotation} (iii), then (i), and \eqref{eq:dmm},
\begin{align}
\Prob(d_{m, m} \geqslant z)
 & = \Prob \biggl( \bigcup_{j=0}^{m-1} \Bigl\{ d^j_{m, m} \geqslant z,\ \min_{l \in \Z} \{ S^j_m(l) - l \} = 0 \Bigr\} \biggr) \nonumber \\
 & \leqslant \sum_{j=0}^{m-1} \Prob\Bigl( d^j_{m, m} \geqslant z,\ \min_{l \in \Z} \{ S^j_m(l) - l \} = 0 \Bigr) \nonumber \\
 & = m \Prob\Bigl( d_{m, m} \geqslant z,\ \min_{l \in \Z} \{ S_m(l) - l \} = 0 \Bigr) \nonumber \\
 & = m \Prob\biggl( \sum_{k=1}^m S_m(k) - k \geqslant z,\ \min_{l \in \Z} \{ S_m(l) - l \} = 0 \biggr) \nonumber \\
 & \leqslant m \Prob\biggl( \sum_{k=1}^m S_m(k) - k \geqslant z \biggr) \label{eq:full_dmm_ub} .
\end{align}

Now, for $\alpha \in \intervalleof{3/2}{2}$ and $\delta \geqslant 0$,
\begin{align}
\Prob\biggl( \sum_{k=1}^m S_m(k) - k \geqslant m^\alpha \delta \biggr)
 & = \Prob\biggl( m \sum_{k=1}^m (L_m - \mathcal{U})(\intervalleff{0}{k/m}) \geqslant m^\alpha \delta \biggr) \nonumber \\
 & = \Prob\biggl( m^2 \biggl( \int_0^1 (L_m - \mathcal U)(\intervalleff{0}{y}) dy + \frac{A_m}{m} \biggr) \geqslant m^\alpha \delta \biggr) \nonumber \\
 & = \Prob\bigl( \varphi(m^{2-\alpha} (L_m - \mathcal U)) \geqslant \delta_m \bigr) \label{eq:full_dmm_phi} ,
\end{align}
where $A_m \in \intervalleff{-1/2}{1/2}$, $\delta_m \defeq \delta + m^{\alpha-1} A_m \to \delta$ almost surely, and, for any measure $\nu \in \mathcal M(\intervalleff{0}{1})$ (the space of signed measures on $\intervalleff{0}{1}$),
\begin{align}
\varphi(\nu)
 \defeq \int_0^1 \nu(\intervalleff{0}{y}) dy
 = \int_0^1 (1 - x) d\nu(x) \label{def:phi}
\end{align}
by Fubini's theorem. In particular, $\varphi$ is a continuous function when $\mathcal M(\intervalleff{0}{1})$ is equipped with the $\tau$-topology, which is generated by the applications $\nu \mapsto \nu(f)$ with $f \colon \intervalleff{0}{1} \to \R$ bounded measurable.

%%%%%%%%%%%%%%%%%%%%%%%%%%%%%%%%%%
\subsection{Upper moderate deviations - Theorem \ref{th:md_full}}
%%%%%%%%%%%%%%%%%%%%%%%%%%%%%%%%%%

Let $\alpha \in \intervalleoo{3/2}{2}$ and $\delta \geqslant 0$. By \cite[Theorem 3.1]{DeAcosta_1994_ProjectiveSystems} (it seems that the result already exists in reference [8] of \cite{Wu94}), we have
\begin{align}
 & - \inf \ensavec{\frac{1}{2} \int_0^1 \Bigl(\frac{d\nu}{dy}(y)\Bigr)^2 dy}{\nu\in \mathcal M(\intervalleff{0}{1}),\ \nu \ll \mathcal U,\ \varphi(\nu) > \delta,\ \nu(\intervalleff{0}{1}) = 0} \label{eq:sanov_md_lb} \\
 & \leqslant \liminf \frac{1}{m^{2\alpha-3}} \log \Prob\bigl(\varphi(m^{2-\alpha}(L_m-\mathcal U)) \geqslant \delta \bigr) \nonumber\\
 & \leqslant \limsup \frac{1}{m^{2\alpha-3}} \log \Prob\bigl(\varphi(m^{2-\alpha}(L_m-\mathcal U)) \geqslant \delta \bigr) \nonumber\\
 & \leqslant - \inf \ensavec{\frac{1}{2} \int_0^1 \Bigl(\frac{d\nu}{dy}(y)\Bigr)^2 dy}{\nu\in \mathcal M(\intervalleff{0}{1}),\ \nu \ll \mathcal U,\ \varphi(\nu) \geqslant \delta,\ \nu(\intervalleff{0}{1}) = 0} \label{eq:sanov_md_up} .
\end{align}

Let us consider the following minimization problem:
\begin{align}
 & \inf \ensavec{\frac{1}{2} \int_0^1 \Bigl(\frac{d\nu}{dy}(y)\Bigr)^2 dy}{\nu\in \mathcal M(\intervalleff{0}{1}),\ \nu \ll \mathcal U,\ \varphi(\nu) = \delta,\ \nu(\intervalleff{0}{1}) = 0}\nonumber\\
 & = \inf \ensavec{\frac{1}{2} \int_0^1 G'(y)^2 dy}{G \in \mathrm{AC}_0(\intervalleff{0}{1}),\ \int_0^1 G(y) dy = \delta} \label{def:pb_var},
\end{align}
where $\mathrm{AC}_0(\intervalleff{0}{1})$ is the space of absolutely continuous functions $G$ on $\intervalleff{0}{1}$ such that $G(0) = G(1) = 0$. Using the method of Lagrange multipliers, if $G$ is a minimizer, then there exists $\lambda \in \R$ such that
\[
\forall h \in \mathrm{AC}_0(\intervalleff{0}{1}) \quad \int_0^1 \bigl( G'(y) h'(y) + \lambda h(y) \bigr) dy = 0 .
\]
Integrating by parts, one has
\[
\forall h \in \mathrm{AC}_0(\intervalleff{0}{1}) \quad \int_0^1 \bigl( G'(y) - \lambda y \bigr) h'(y) dy = 0 .
\]
By Du Bois-Reymond's lemma in \cite[p.184]{Clarke13}, the function $y \mapsto G'(y) - \lambda y$ is constant, so $G$ is a quadratic polynomial. Getting back to \eqref{def:pb_var}, the minimizer is given by
\[
G(y) = 6\delta y(1-y)
\]
and the infimum is $6 \delta^2$. Consequently, using the continuity of $\delta \mapsto 6 \delta^2$ to lower bound \eqref{eq:sanov_md_lb} and the positive homogeneity of the constraints in \eqref{eq:sanov_md_up}, we get
\begin{align} \label{eq:full_md_phi_limit}
\frac{1}{m^{2\alpha-3}} \log \Prob\bigl(\varphi(m^{2-\alpha}(L_m-\mathcal U)) \geqslant \delta \bigr) \to - 6 \delta^2.
\end{align}

Using the fact that, for $\Phi = \varphi(m^{2-\alpha}(L_m-\mathcal U))$, 
\begin{align}\label{eq:argument_cont}
\Prob(\Phi\geqslant \delta+\varepsilon)
&\leqslant 
\Prob(\Phi\geqslant \delta_m)
\leqslant \Prob(\Phi\geqslant \delta-\varepsilon)
\end{align}
for all $\varepsilon>0$ and all $m$ large enough, Theorem \ref{th:md_full} stems from \eqref{eq:full_dmm_lb}, \eqref{eq:full_dmm_ub}, \eqref{eq:full_dmm_phi}, \eqref{eq:full_md_phi_limit}, \eqref{eq:argument_cont}, and by the continuity of the function $\delta \mapsto -6\delta^2$.\qed
%\end{proof}

%%%%%%%%%%%%%%%%%%%%%%%%%%%%%%%%%%
\subsection{Upper large deviations - Theorem \ref{th:ld_full}}
%%%%%%%%%%%%%%%%%%%%%%%%%%%%%%%%%%

%\alert{Regarder si on peut montrer la borne supérieure avec Spencer 2000, ce qui simplifie l'écriture de notre borne inférieure.}

%\begin{proof}[Proof of Theorem \ref{th:ld_full}]\leavevmode 
Here $\alpha = 2$. Let $\delta \geqslant 0$. The result for $\delta = 0$ is trivial. It is also trivial for $\delta \geqslant 1/2$, since $d_{m, m} \leqslant m(m-1)/2 < m^2/2$. Assume that $\delta \in \intervalleoo{0}{1/2}$. By Sanov's theorem (see, e.g., \cite{Csiszar_2006_Sanov}), we have
\begin{align}
 & - \inf \ensavec{ \int_0^1 \Bigl(\frac{d\nu}{dy}(y)\Bigr)\log\Bigl(\frac{d\nu}{dy}(y)\Bigr)dy}{\nu\in \mathcal M_1^+(\intervalleff{0}{1}),\ \nu \ll \mathcal U,\ \varphi(\nu) > \delta + 1/2} \label{eq:sanov_ld_lb} \\
 & \leqslant \liminf \frac{1}{m} \log \Prob\left(\varphi\left(L_m-\mathcal U\right)\geqslant \delta\right) \nonumber\\
 & \leqslant \limsup \frac{1}{m} \log \Prob\left(\varphi\left(L_m-\mathcal U\right)\geqslant \delta\right)\nonumber\\
 & \leqslant - \inf \ensavec{ \int_0^1 \Bigl(\frac{d\nu}{dy}(y)\Bigr)\log\Bigl(\frac{d\nu}{dy}(y)\Bigr)dy}{\nu\in \mathcal M_1^+(\intervalleff{0}{1}),\ \nu \ll \mathcal U,\ \varphi(\nu) \geqslant \delta + 1/2} \label{eq:sanov_ld_up} ,
\end{align}
where $\mathcal M_1^+(\intervalleff{0}{1})$ is the space of probability measures on $\intervalleff{0}{1}$. Let us consider the following minimization problem:
\begin{align}
 & \inf\ensavec{ \int_0^1 \Bigl(\frac{d\nu}{dy}(y)\Bigr)\log\Bigl(\frac{d\nu}{dy}(y)\Bigr)dy}{\nu\in \mathcal M_1^+(\intervalleff{0}{1}),\ \nu \ll \mathcal U,\ \varphi(\nu) = \delta + 1/2} \nonumber\\
 & = \inf\biggl\{ \int_0^1 F'(y)\log F'(y)dy\ ; \ F\in \mathrm{AC}(\intervalleff{0}{1}),\ F'\geqslant 0,\ \int_0^1 F(y) dy = \delta + 1/2,\nonumber\\
 & \qquad \qquad \qquad  F(0)=0,\ F(1)=1\biggr\} \nonumber \\
 & = \inf K \label{def:pb_var_2} ,
\end{align}
where $\mathrm{AC}(\intervalleff{0}{1})$ is the space of absolutely continuous functions on $\intervalleff{0}{1}$ and $K \colon \mathrm{AC}(\intervalleff{0}{1}) \to \intervalleff{0}{\infty}$ is the convex function defined by
\begin{align*}
K(F)
 = \begin{cases}
\int_0^1 F'(y)\log(F'(y))dy & \text{if $F'\geqslant 0$, $\int_0^1 F(y)dy = \delta+\frac 12$, $F(0)=0$, $F(1)=1$} \\
\infty & \text{otherwise.}
\end{cases}
\end{align*}
It is a standard convex optimization problem, a minimizer of which is
\[
\bar F(y)=a(1-e^{\lambda y}), \quad  \text{where} \quad 
\begin{cases} 
a(1-e^{\lambda})=1\\
a=\delta+\frac 12-\frac 1\lambda.
\end{cases}
\]
(One can see that $a > 1$ and $\lambda < 0$.) Indeed, by the definition of a convex function and the subdifferential, it suffices to check that $0$ belongs to the subdifferential of $K$ at $\bar F$. For all $h\in \mathrm{AC}(\intervalleff{0}{1})$ and for all $t>0$, 
\[
K(\bar F+th)=\begin{cases}
\int_0^1 (\bar F'+th')\log(\bar F'+th') & \text{if $\bar{F}'+th'\geqslant 0$, $\int_0^1 h = 0$, $h(0)=h(1)=0$} \\
\infty& \text{otherwise.}
\end{cases}
\]
Differentiating under the integral sign with respect to $t$ and integrating by parts gives
\[
K'(\bar F;h)
 = \int_0^1 h'(y)(\log(\bar{F}'(y)) + 1)dy
 = -\lambda \int_0^1 h(y) dy
 = 0 ,
\]
since $h(0) = h(1) = 0$ and $\int_0^1 h(y) dy = 0$. It remains to compute the value of $K$ at $\bar F$ and to conclude following the same arguments (continuity and positive homogeneity of the constraints) as in the end of the proof of Theorem \ref{th:md_full}.\qed 
%\end{proof}

\section{Interlude}\label{sec:inter}

\subsection{Janson's reformulation}

%Notation $X$, $Y$, $S_n$ et $T_n$ avec des $i$ pour les sommes\\
%Lien entre $d_{m,n}$ et les va ci-dessus. Cf. Janson.\\
%Donner les deux définitions de la loi de Borel et dire que $\mu_n$ et $\lambda_n$ convergent vers ...

%Notice that $m = \Theta(n)$ and $N = m-n = \Theta(m)$.

Here, we consider $(m,n)$ with $m \to \infty$ and $n/m \to \mu \in \intervalleoo{0}{1}$. As a consequence, $N = m-n \to \infty$. To make notation clearer, we make quantities depend on $N$ and all limits are considered as $N\to \infty$ unless stated otherwise.
%Let $(m_n)_{n \geqslant 1}$ be a sequence of positive integers. Assume that $n/m_n \to \mu \in \intervalleoo{0}{1}$.
In the next section, we are interested in the deviations of $d_{m_N,n_N}$ in that regime, which is called the sparse case (see \cite{FPV98,Janson01a} for this denomination with slight variants). In this section, we introduce a reformulation of the model of hashing with linear probing due to Janson in \cite{Janson01a} and prove some preliminary results.

\medskip
%\alert{Commencer ici à changer $n$ en $N$ (et $N_n$ en $N$).}

For all $N \geqslant 1$, we consider a vector of random variables $(\rva{X}{N},\rva{Y}{N})$ defined as follows. We assume that $\rva{X}{N}$ is distributed according to the Borel distribution with parameter $\mu_N \defeq n_N/m_N \in \intervalleoo{0}{1}$, i.e.
\begin{align}\label{def:X_mu}
\forall l \in \intervallentfo{1}{\infty} \quad \Prob(\rva{X}{N} = l) = e^{-\mu_N l} \frac{(\mu_N l)^{l-1}}{l!}
\end{align}
(see, e.g., \cite{FPV98} or \cite{Janson01a} for more details). 
In some places, for the ease of computation, we may also use the parametrization $\lambda_N=e^{-\mu_N} \mu_N$ to get an equivalent definition of the Borel distribution:
\begin{align}\label{def:X_lambda}
\Prob(\rva{X}{N}=l)=\frac{1}{T(\lambda_N)}\frac{l^{l−1}\lambda_N^l}{l!},
\end{align}
where $T$ is the tree function  (see, e.g., \cite[p.\ 127]{FS09}). 
%Naturally, the convergence of $\mu_N$ to $\mu$ implies the one of $\lambda_N$ to $\lambda\defeq  e^{-\mu} \mu$. 
Furthermore, we assume that $\rva{Y}{N}$ given $\{ \rva{X}{N} = l \}$ is distributed as $d_{l,l-1}$.

%Now, let $N \defeq m_N - n$ (notice that $N \to \infty$). 
Let $(\rva[i]{X}{N},\rva[i]{Y}{N})_{ 1\leqslant i\leqslant N}$ be an i.i.d.\ sample distributed as $(\rva{X}{N},\rva{Y}{N})$ and define, for all $k \in \intervallentff{1}{N}$,
\[
\rva[k]{S}{N} \defeq \sum_{i=1}^{k} \rva[i]{X}{N}
\quad \text{and} \quad
\rva[k]{T}{N}\defeq\sum_{i=1}^{k} \rva[i]{Y}{N} .
\]
To lighten notation, let $\rva{S}{N} \defeq \rva[N]{S}{N}$ and $\rva{T}{N} \defeq \rva[N]{T}{N}$. Notice that, for all $N \geqslant 1$, $\Espe[\rva{X}{N}]=(1-\mu_N)^{-1}=m_N/N$, so $\Espe[\rva{S}{N}]=m_N$. Moreover, $\Prob(\rva{S}{N}= m_N) > 0$ and we have the following identity (see \cite[Lemma 4.1]{Janson01a}):
\[
\mathcal{L}(d_{m_N,n_N}) = \mathcal{L}( \rva{T}{N} \, | \, \rva{S}{N} = m_N) .
\]

\subsection{Tail estimates}

For all $\xi\in \intervalleoo{0}{1}$, recall that
$\kappa(\xi)=\xi - \log(\xi) - 1 \in \intervalleoo{0}{\infty}$.

\begin{prop}[Tail of $\rva{X}{N}$] \label{prop:queue_X}
If $l \geqslant 1/\mu_N$, then
\begin{align}\label{eq:queue_X_1}
\log \Prob(\rva{X}{N} = l) \leqslant - \kappa(\mu_N) l .
\end{align}
And if $l_N \to \infty$, then
\begin{align}\label{eq:queue_X_2}
\log \Prob(\rva{X}{N} \geqslant l_N) \sim \log \Prob(\rva{X}{N} = l_N) \sim -\kappa(\mu) l_N .
\end{align}
\end{prop}

\begin{proof}[Proof of Proposition \ref{prop:queue_X}]
As soon as $\mu_N l \geqslant 1$, and since $\log(l!) \geqslant l(\log(l)-1)$,
\begin{align*}
\log \Prob(\rva{X}{N}=l)& = -\mu_N l + l\log(\mu_N l) - \log(\mu_N l) - \log(l!) 
 \leqslant -l(\mu_N-\log(\mu_N)-1)
 = - \kappa(\mu_N) l .
\end{align*}
Therefore, 
\begin{align*}
\log \Prob(\rva{X}{N}\geqslant l_N)
= \log \sum_{l=l_N}^{\infty} \Prob(\rva{X}{N}=l_N)
\leqslant \log \sum_{l=l_N}^{\infty} e^{ - \kappa(\mu_N) l}\sim - \kappa (\mu)l_N.
\end{align*}
Finally, using Stirling formula, one has
\begin{align*}
\log \Prob(\rva{X}{N}\geqslant l_N)\geqslant \log \Prob(\rva{X}{N}=l_N) = -\mu_N l_N + (l_N-1)\log(\mu_N l_N) - \log(l_N!) \sim - \kappa(\mu) l_N. \quad \qedhere
\end{align*}
\end{proof}

From the previous proposition and Theorem \ref{th:ld_full}, we deduce the asymptotic behavior of the tail of the pair $(\rva{X}{N},\rva{Y}{N})$.

\begin{prop}[Tail of $(\rva{X}{N},\rva{Y}{N})$] \label{prop:queue_couple}
Let $l_N \to \infty$ and let $p_N$ be such that $p_N/l_N^2 \to \delta$. Then
\begin{align}\label{eq:queue_couple}
\frac{1}{l_N} \log \Prob(\rva{X}{N}=l_N,\ \rva{Y}{N} \geqslant p_N) \to -(\kappa(\mu) + J(\delta)) .
\end{align}
\end{prop}

\begin{proof}[Proof of Proposition \ref{prop:queue_couple}]
It suffices to write
\begin{align*}
\frac{1}{l_N} \log \Prob(\rva{X}{N}=l_N,\ \rva{Y}{N} \geqslant p_N)
 & = \frac{1}{l_N} \log  \Prob(\rva{X}{N}=l_N) + \frac{1}{l_N} \log \Prob(d_{l_N, l_N-1} \geqslant p_N) \\
 & = \frac{1}{l_N} \log  \Prob(\rva{X}{N}=l_N) + \frac{1}{l_N} \log \Prob(d_{l_N, l_N} \geqslant l_N^2(\delta + o(1))) \\
 & \to - (\kappa(\mu) + J(\delta))
\end{align*}
by Remark \ref{rem:dense_case}, Proposition \ref{prop:queue_X}, and Theorem \ref{th:ld_full}.
\end{proof}

\begin{lem} \label{lem:maj_queue_XY_uniforme}
Let $\tilde{J} \leqslant J$ be any nondecreasing function, continuous on $\intervalleff{0}{1/2}$. For all $\varepsilon > 0$, there exists $N_0 \geqslant 1$ and $l_0 \geqslant 1$ such that, for all $N \geqslant N_0$, for all $l \geqslant l_0$, and for all $\delta \geqslant 0$,
\[
\log \Prob(\rva{X}{N}=l,\ \rva{Y}{N} \geqslant \delta l^2) \leqslant -(\kappa(\mu) + \tilde{J}(\delta)-\varepsilon)l.
\]
\end{lem}

\begin{proof}[Proof of Lemma \ref{lem:maj_queue_XY_uniforme}]
The result is trivial for $\delta \in \intervalleoo{1/2}{\infty}$. Remember that
\[
\Prob(\rva{X}{N}=l,\ \rva{Y}{N} \geqslant \delta l^2) =  \Prob(\rva{X}{N}=l) \Prob(d_{l,l-1} \geqslant \delta l^2) .
\]
On the one hand, let $\varepsilon > 0$. By Proposition \ref{prop:queue_X}, if $N$ and $l$ are large enough,
\begin{align*}
\log \Prob(\rva{X}{N}=l)
 \leqslant - \kappa(\mu_N) l
 \leqslant -(\kappa(\mu)-\varepsilon/2)l .
\end{align*}
On the other hand, the nondecreasing functions 
$\phi_l\colon \delta \in \intervalleff{0}{1/2} \mapsto \min(-l^{-1}\log\Prob(d_{l,l-1} \geqslant \delta l^2), \tilde{J}(\delta))$ converge pointwise to  the continuous function $\tilde{J}|_{\intervalleff{0}{1/2}}$ as $l \to \infty$, by Theorem \ref{th:ld_full} ; thus the convergence is uniform and the result follows.
\end{proof}

In the sequel, we will also need the asymptotic behavior of the tail of $\rva{Y}{N}$ alone.

\begin{prop}[Tail of $\rva{Y}{N}$] \label{prop:queue_Y}
If $p_N \to \infty$, then
\begin{align}\label{eq:queue_Y}
\frac{1}{\sqrt{p_N}} \log \Prob(\rva{Y}{N} \geqslant p_N) \to -q(\mu) = -\inf_{0<\delta<1/2} \frac{1}{\sqrt{\delta}}\bigl[\kappa(\mu) + J(\delta)\bigr].
\end{align}
\end{prop}

\begin{proof}[Proof of Proposition \ref{prop:queue_Y}]
For $\delta > 0$, let $l_N = \ceil{(p_N/\delta)^{1/2}}$. Then,
\begin{align*}
\frac{1}{\sqrt{p_N}} \log \Prob(\rva{Y}{N} \geqslant p_N)
 \geqslant \frac{l_N}{\sqrt{p_N}} \cdot \frac{1}{l_N} \log\Prob(\rva{X}{N} = l_N,\ \rva{Y}{N} \geqslant p_N)
  \to  -\frac{1}{\sqrt{\delta}}\bigl[\kappa(\mu)+J(\delta)\bigr] ,
\end{align*}
by Proposition \ref{prop:queue_couple}. Taking the supremum in $\delta > 0$, one gets
\begin{equation} \label{eq:queue_Y_minoration}
\liminf_{N \to \infty} \frac{1}{\sqrt{p_N}} \log \Prob(\rva{Y}{N} \geqslant p_N) \geqslant - \inf_{0 < \delta < 1/2} \frac{1}{\sqrt{\delta}} \bigl[\kappa(\mu) + J(\delta)\bigr] .
\end{equation}

Now we turn to the upper bound. Let us fix $\beta > 0$ such that $\beta \kappa(\mu) > q(\mu)$. Let $l_N = \floor{\beta p_N^{1/2}}$ and write
\begin{align*}
\Prob(\rva{Y}{N} \geqslant p_N)
 = \sum_{l=1}^{l_N} \Prob(\rva{X}{N} = l,\ \rva{Y}{N} \geqslant p_N) +  \sum_{l=l_N+1}^\infty \Prob(\rva{X}{N} = l,\ \rva{Y}{N} \geqslant p_N)
 \eqdef P_N + R_N .
\end{align*}
First of all, using Proposition \ref{prop:queue_X},
\[
\frac{1}{\sqrt{p_N}} \log(R_N)
 \leqslant \frac{1}{\sqrt{p_N}} \log\Prob(\rva{X}{N} > l_N)
  \to  - \beta \kappa(\mu) < - q(\mu) .
\]
Let $\varepsilon > 0$. Taking into account the already proved lower bound, and using Lemma \ref{lem:maj_queue_XY_uniforme} with
\[
\tilde{J}(\delta) = J_\varepsilon(\delta)
 = \begin{cases}
J(\delta) \wedge \varepsilon^{-1} & \text{if $\delta \leqslant 1/2$} \\
\infty & \text{if $\delta > 1/2$,}
\end{cases}
\]
we deduce that
\begin{align*}
\limsup_{N \to \infty} \frac{1}{\sqrt{p_N}} \log \Prob(\rva{Y}{N} \geqslant p_N)
 & = \limsup_{N \to \infty} \frac{1}{\sqrt{p_N}} \log(P_N) \\
 & \leqslant \max_{1 \leqslant l \leqslant l_N} - \frac{l}{\sqrt{p_N}} \Bigl[ \kappa(\mu) + J_\varepsilon\Bigl(\frac{p_N}{l^2}\Bigr) - \varepsilon \Bigr] \\
 & \leqslant - \inf_{1/\beta \leqslant \delta < 1/2} \frac{1}{\sqrt{\delta}} \bigl[ \kappa(\mu) + J_\varepsilon(\delta) - \varepsilon \bigr] \\
 & \eqdef M_\varepsilon .
\end{align*}

Let $\delta_\varepsilon \in \intervalleoo{0}{1/2}$ be such that $J(\delta_\varepsilon) = 1/\varepsilon$. We have
\begin{align*}
\inf_{\delta_\varepsilon < \delta < 1/2} & \frac{1}{\sqrt{\delta}} \bigl[ \kappa(\mu) + J_\varepsilon(\delta) - \varepsilon \bigr]
 \geqslant \sqrt{2} \bigl( \kappa(\mu) + \varepsilon^{-1} - \varepsilon \bigr)
 \xrightarrow[\varepsilon \to 0]{} \infty .
\end{align*}
A fortiori, since $J_\varepsilon \leqslant J$,
\begin{align*}
\inf_{\delta_\varepsilon < \delta < 1/2} & \frac{1}{\sqrt{\delta}} \bigl[ \kappa(\mu) + J(\delta) - \varepsilon \bigr]
 \xrightarrow[\varepsilon \to 0]{} \infty .
\end{align*}
So, if $\varepsilon$ is small enough,
\begin{align*}
M_\varepsilon
& = - \inf_{1/\beta \leqslant \delta < 1/2} \frac{1}{\sqrt{\delta}} \bigl[ \kappa(\mu) + J(\delta) - \varepsilon \bigr]\\
& \xrightarrow[\varepsilon \to 0]{} - \inf_{1/\beta \leqslant \delta < 1/2} \frac{1}{\sqrt{\delta}} \bigl[ \kappa(\mu) + J(\delta) \bigr] \leqslant - \inf_{0< \delta < 1/2} \frac{1}{\sqrt{\delta}} \bigl[ \kappa(\mu) + J(\delta) \bigr]
\end{align*}
and the result follows.
\end{proof}

\subsection{Useful limit theorems}\label{sec:inter_res}

The following lemma is a direct consequence of \cite[Lemma 4.3]{Janson01a} and Proposition \ref{prop:cv_moments}.
 
\begin{prop}\label{prop:moy_weak_array_cond_mob_cor_moderate} 
One has
\[
\Espe[ \rva{T}{N} \, | \,\rva{S}{N}= m_N] = \Espe[\rva{T}{N}] + o(N^{1/2}).
\]
\end{prop}

Let $(X,Y)$ be a pair of random variables such that $X$ is distributed according to the Borel distribution with parameter $\mu=\lim \mu_N$ and $Y$ given $\{ X = l \}$ is distributed as $d_{l,l-1}$. Let $\lambda=e^{-\mu}\mu$ be the other standard parameter of the Borel distribution as in \eqref{def:X_lambda}.

\begin{prop}[Moments convergence] \label{prop:cv_moments}
$(\rva{X}{N},\rva{Y}{N})_{N \geqslant 1}$ converges to  $(X,Y)$ in distribution and with all mixed moments of the type $\Espe[\rva{X}{N}^p \rva{Y}{N}^qe^{(s+it)\rva{X}{N}}]$, where $p\geqslant 0$, $q\geqslant 0$,  $s<-\log(\lambda e)$, and $t\in \R$.
\end{prop}

%\begin{remark}
%Cf. article Janson \cite[Remark 1.2]{Janson01a}.
%\end{remark}

\begin{proof}[Proof of Proposition \ref{prop:cv_moments}]
Let $f\colon \R\times \R \to \R$ be a bounded measurable function. Using \eqref{def:X_lambda}, one has 
\begin{align*}
\Espe[f(\rva{X}{N},\rva{Y}{N})]&=\sum_{l \in \N} \Espe[f(\rva{X}{N},\rva{Y}{N})\, | \, \rva{X}{N} = l] \Prob(\rva{X}{N} = l)\\
&=\sum_{l \in \N} \Espe[f(l,d_{l,l-1})] \frac{1}{T(\lambda_N)} \frac{l^{l-1}}{l!} \lambda_N^l.
\end{align*}
Since $\lambda_N$ converges to $\lambda$, $T$ is continuous, and
\begin{align}\label{eq:tree}
\frac{1}{T(\lambda_N)} \frac{l^{l-1}}{l!} \lambda_N^l  \leqslant \frac{T(\lambda+\varepsilon)}{T(\lambda-\varepsilon)} \biggl( \frac{1}{T(\lambda+\varepsilon)} \frac{l^{l-1}}{l!} (\lambda+\varepsilon)^l \biggr),
\end{align}
as soon as $\abs{\lambda_N - \lambda} \leqslant \varepsilon$, we conclude by Lebesgue's dominated convergence theorem that    
$(\rva{X}{N},\rva{Y}{N})_{N \geqslant 1}$ converges in distribution to $(X,Y)$, where $X$ is Borel distributed with parameter $\lambda$ and $\mathcal L(Y\, | \,X=l)=\mathcal L(d_{l,l-1})$. 

Let $c>1$ such that $sc<-\log(\lambda e)$ and $(a,b)\in (\R_+^*)^2$ such that $a^{-1}+b^{-1}+c^{-1}=1$. Hölder's inequality yields
\begin{align*}
\lvert \Espe[\rva{X}{N}^p\rva{Y}{N}^qe^{(s+it)\rva{X}{N}}]\rvert\leqslant
\Espe[\rva{X}{N}^{ap}]^{1/a}\Espe[\rva{Y}{N}^{bq}]^{1/b}\Espe[e^{sc\rva{X}{N}}]^{1/c}.
\end{align*}
By \eqref{eq:tree}, for each $r > 0$ and $s'<-\log(\lambda e)$, $\limsup \Espe[\rva{X}{N}^r]$ and $\limsup \Espe[e^{s'\rva{X}{N}}]$ are finite.  
Moreover, since $d_{ l,l-1} \leqslant l^2$,
\begin{align*}
\Espe[\rva{Y}{N}^r]
 & = \sum_{l \in \N} \Espe[\rva{Y}{N}^r\, | \, \rva{X}{N} = l] \Prob(\rva{X}{N} = l) 
  \leqslant \sum_{l \in \N} l^{2r} \Prob(\rva{X}{N} = l) 
  = \Espe[\rva{X}{N}^{2r}] ,
\end{align*}
so $\limsup \Espe[\rva{Y}{N}^r]$ is finite too. Hence, by uniform integrability (see, e.g., \cite[Example 2.21]{van1998asymptotic}), we obtain the convergence of all mixed moments.
\end{proof}

\begin{prop}[Local large deviations for $\rva{S}{N}$]\label{prop:LLT_gd_dev}
For any sequence of integers$(k_N)_{N \geqslant 1}$ such that $\lim k_N/N\in \intervalleoo{1}{\infty}$, we have
\[
\log \Prob(\rva{S}{N} = k_N) = - N \Lambda_{\rva{X}{N}}^*(k_N/N) + O(\log(N)).
\]
\end{prop}

%\alert{REMARQUE : on utilise :
%\begin{itemize}
%\item $\Lambda_{\rva{X}{N}}^*(m_N/N) = 0$ ;
%\item $\Lambda_{\rva{X}{N}}^*((m_N - l_N)/N) = \frac{1}{2 \sigma_X^2} (l_N ∕ N)^2 + O((l_N ∕ N)^3)$.
%\end{itemize}}

\begin{proof}[Proof of Proposition \ref{prop:LLT_gd_dev}]
We just check that we can apply \cite[Lemma 3.3]{TFC12} to the sequence $(\rva{X}{N})_{N \geqslant 1}$. The conclusion follows since, in this case, $m = 1$, $b = 0$ and $c_{n, m, b} = 1$.
First, $\Ima(\Lambda_{\rva{X}{N}}')=\intervalleoo{1}{\infty}=\Ima(\Lambda_X')$ so, for all $N$ large enough, $k_N/N\in \Ima(\Lambda_{\rva{X}{N}}')$ and $\lim k_N/N\in \Ima(\Lambda_X')$. Secondly, $\interieur(\dom(\Lambda_{\rva{X}{N}}))=\intervalleoo{-\infty}{-\log(\lambda_N e)}$ and $\interieur(\dom(\Lambda_X))=\intervalleoo{-\infty}{-\log(\lambda e)}$, so that assumption 1.\ of \cite[Lemma 3.3]{TFC12} holds for all $N$ large enough (since $\lambda_N\to \lambda$). Thirdly, assumption 2.\ of \cite[Lemma 3.3]{TFC12} stems from Proposition \ref{prop:cv_moments}. 
\end{proof}

The following proposition is a non conditioned version of Theorem \ref {th:umld_sparse}. It stems immediately from \cite{ATP2020array} (with $\epsilon=1/2$ and $q=q(\mu)$, defined in Theorem \ref{th:umld_sparse}, (ii)) and Propositions \ref{prop:queue_Y} and \ref{prop:cv_moments}.

\begin{prop}[Large deviations for $\rva{T}{N}$]\label{prop:lognag}
\leavevmode
\begin{enumerate}
\item[(i)]
If $\alpha<2/3$, then
\begin{align}\label{eq:lognag_dev_mod}
\lim_{N \to \infty} \frac{1}{N^{2\alpha-1} }\log \Prob( \rva{T}{N}-\Espe[\rva{T}{N}] \geqslant N^{\alpha} y ) = - \frac{y^2}{2\sigma^2(\mu)}.
\end{align}

\item[(ii)] 
If $\alpha=2/3$, then
\begin{align}\label{eq:lognag_dev_inter}
\lim_{N \to \infty} \frac{1}{N^{1/3} }\log \Prob( \rva{T}{N} -\Espe[\rva{T}{N}] \geqslant N^{2/3} y ) = - I(y)
\end{align}
where $I$ is defined in Theorem \ref{th:umld_sparse}, (ii).
%\[
%I(y)\defeq \begin{cases}
%\frac{y^2}{2\sigma^2(\mu)} & \text{if  $y\leqslant y(\mu)$},\\ 
%q(\mu)(1-t(y))^{1/2} y^{1/2}+\frac{t(y)^2y^2}{2\sigma^2(\mu)} & \text{if  $y> y(\mu)$}, 
%\end{cases}
%\]
%with $y(\mu)=3\left(q(\mu)\sigma^2(\mu)\right)^{2/3}/2$ and  $t(y)$ being defined for $y>y(\mu)$ as the smallest root of the cubic equation in $t\in \intervalleff{0}{1}$: 
%\[
%t^3-t^2+\frac{q^2(\mu)\sigma^4(\mu)}{4y^{3}}=0.
%\]

\item[(iii)]If  $\alpha>2/3$, then
\begin{align}\label{eq:lognag_gdes_dev}
\lim_{N \to \infty} \frac{1}{N^{\alpha/2} }\log \Prob( \rva{T}{N}-\Espe[\rva{T}{N}]  \geqslant N^{\alpha} y ) = - q(\mu) y^{1/2}.
\end{align}
\end{enumerate}
\end{prop}

%%%%%%%%%%%%%%%%%%%%%%%%%%%%%%%%%%
%%%%%%%%%%%%%%%%%%%%%%%%%%%%%%%%%%
%%%%%%%%%%%%%%%%%%%%%%%%%%%%%%%%%%
\section{Proofs: sparse tables}\label{sec:proof_sparse}
%%%%%%%%%%%%%%%%%%%%%%%%%%%%%%%%%%
%%%%%%%%%%%%%%%%%%%%%%%%%%%%%%%%%%
%%%%%%%%%%%%%%%%%%%%%%%%%%%%%%%%%%

%%%%%%%%%%%%%%%%%%%%%%%%%%%%%%%%%%
\subsection{Lower moderate deviations - Theorem \ref{th:lmd_sparse}}
%%%%%%%%%%%%%%%%%%%%%%%%%%%%%%%%%%

%\begin{proof}[Proof of ]
%Ingredients :
%\begin{itemize}
%\item reformulation en $(X_N, Y_N)$ ;
%\item Plachky-Steinbach ;
%\item calculs de \cite{TFC12}.
%\end{itemize}
One has
\begin{align}
\Prob(d_{m_N,n_N} - \Espe[d_{m_N,n_N}] \leqslant - N^\alpha y)
 & = \Prob( T_N - \Espe[T_N \, | \, S_N=m_N] \leqslant - N^\alpha y \, |\, S_N=m_N ) \nonumber\\
 & = \Prob(T_N - \Espe[T_N] \leqslant - N^\alpha y_N \, |\, S_N=m_N)\label{eq:prob}
\end{align}
where
\[
y_N \defeq  y - \frac{1}{N^\alpha}( \Espe[T_N \, | \, S_N=m_N] - \Espe[T_N] )  \to  y
\]
by Proposition \ref{prop:moy_weak_array_cond_mob_cor_moderate}. Since the variables are nonnegative, their Laplace transforms are defined on $\intervalleoo{-\infty}{0}$ at least. Adapting the proof of \cite[Theorem 2.2]{TFC12} to the unilateral case and using \cite{PS75} (unilateral version of Gärtner-Ellis theorem), we get \eqref{eq:lmd_sparse}. \qed
%\[
%\pushQED{\qed} 
%\frac{1}{N^{2\alpha-1}} \log \Prob(d_{m_N,n_N} - \Espe[d_{m_N,n_N}] \leqslant - N^{\alpha} y)
%\to  - \frac{y^2}{2\sigma^2(\mu)}.\qedhere
%\popQED
%\]
%Éventuellement à préciser :
%\begin{itemize}
%\item $m_N/N = 1/(1-\mu_N) \in R_{X^{(N)}} = \intervalleoo{0}{\infty}$ ;
%\item $1/(1-\mu) \in R_{X} = \intervalleoo{0}{\infty}$ ;
%\item $\ensavec{f'(u)}{u \leqslant 0} \supset \intervallefo{0}{\infty} \supset \intervallefo{0}{\Espe[Y]}$.
%\end{itemize}
%\end{proof}

%%%%%%%%%%%%%%%%%%%%%%%%%%%%%%%%%%
\subsection{Lower large deviations - Theorem \ref{th:lld_sparse}}
%%%%%%%%%%%%%%%%%%%%%%%%%%%%%%%%%%

For any $\R^d$-valued random variable $Z$, we denote by $\Lambda_Z$ the log-Laplace transform of $Z$, i.e.\ the function defined,  for $\lambda \in \R^d$, by 
\[
\Lambda_Z(\lambda)=\log \Espe[\exp(\lambda\cdot{} Z)],
\]
and by $\Lambda_Z^*$ the Fenchel-Legendre transform of the function $\Lambda_Z$, i.e.\ the function defined, for $z \in \R^d$, by 
\[
\Lambda_Z^*(z)=\sup \ensavec{\lambda\cdot{} z-\Lambda_Z(\lambda)}{\lambda \in \R^d}.
\]

%\begin{proof}[Proof of Theorem \ref{th:lld_sparse}]
%Ingredients :
%\begin{itemize}
%\item reformulation en $(X_N, Y_N)$ ;
%\item Plachky-Steinbach ;
%\item calculs de \cite{TFC12}.
%\end{itemize}
Proceeding as in the proof of Theorem \ref{th:lmd_sparse}, we get 
%\eqref{eq:prob} with $\alpha=1$, one has
%\begin{align*}
%\Prob(d_{m_N,n_N} - \Espe[d_{m_N,n_N}] \leqslant - N y)
% & = \Prob(T_N - \Espe[T_N] \leqslant - N y_N \, |\, S_N=m_N).
%\end{align*}
%Since the variables are nonnegative, their Laplace transforms are defined on $\intervalleoo{-\infty}{0}$. 
%Adapting the proof of \cite[Theorem 2.1]{TFC12} to the unilateral case, and thus using \cite{PS75} instead of Gärtner-Ellis theorem, we get:
\begin{align}
\frac{1}{N} \log \Prob(d_{m_N,n_N} - \Espe[d_{m_N,n_N}] \leqslant - N y)
 &  \to  - \Lambda_{(X,Y)}^*\Bigl( \frac{1}{1-\mu}, \frac{\mu^2}{2(1-\mu)^2}-y \Bigr), %\eqdef - K(y),
\end{align}
since $m_N/N \to 1/(1-\mu)$ and admitting that $\Lambda_{(X,Y)}^*((1-\mu)^{-1},\cdot{})$ is strictly convex on $\intervalleoo{0}{\mu^2/(2(1-\mu)^2)}$. 

Let $x_0=(1-\mu)^{-1}$. Let us prove that $\Lambda_{(X,Y)}^*(x_0,\cdot{})$ is strictly convex. Let $y\in \intervalleoo{0}{\mu^2/(2(1-\mu)^2)}$. First, we prove that $\Lambda^*_{(X,Y)}$ is differentiable at $(x_0,y)$.  
We have $(x_0,y)\in \interieur(\dom(\Lambda_{(X,Y)}^*))$, therefore the subdifferential $\partial \Lambda_{(X,Y)}^*(x_0,y)$ is nonempty, i.e.\ there exists $(s,t) \in \partial \Lambda_{(X,Y)}^*(x_0,y)$ (see \cite[Theorem 23.4]{Rockafellar70convex}). It remains to prove that such a point $(s,t)$ is unique.
Choosing $\varepsilon>0$ such that $\Lambda_{Y}^*(y+\varepsilon)>0$,  
\begin{align*}
-\Lambda_{(X,Y)}^*(x_0,y)
&\leqslant \liminf \frac{1}{N} \log \Prob(\rva{S}{N}\leqslant N(x_0+\varepsilon),\rva{T}{N}\leqslant N(y+\varepsilon) ) \\
&\leqslant \liminf \frac{1}{N} \log \Prob(\rva{T}{N}\leqslant N(y+\varepsilon) ) = -\Lambda_{Y}^*(y+\varepsilon)<0.
\end{align*}
Since $\Lambda_{(X,Y)}^*(x_0,\mu^2/(2(1-\mu)^2))=0$ and $\Lambda^*(x_0,\cdot{})$ is convex, one has $t<0$. Therefore $(s,t)\in \interieur(\dom(\Lambda_{(X,Y)}))$. 
To obtain a local version of \cite[Theorem 23.5]{Rockafellar70convex}, we notice that
\[
(s,t) \in \partial \Lambda_{(X,Y)}^*(x_0,y) \iff (x_0,y) \in \partial \Lambda^{**}_{(X,Y)}(s,t)=\partial \Lambda_{(X,Y)}(s,t)=\{\nabla \Lambda_{(X,Y)}(s,t)\},
\]
since $\Lambda_{(X,Y)}$ is differentiable on $\interieur(\dom(\Lambda_{(X,Y)}))$. 
Now,
\[
\det (\Hess(\Lambda_{(X,Y)})(\lambda,\rho))=\Var (\tilde X)\Var(\tilde Y)-\Cov(\tilde X,\tilde Y)^2>0
\]
where $(\tilde X,\tilde Y)$ has a mass function proportional to $e^{\lambda x+\rho y}f_{(X,Y)}(x,y)$ which is not supported by a line, so $\Lambda_{(X,Y)}$ is strictly convex. Thus $(s,t)$ is the unique solution of $(x_0,y)=\nabla \Lambda_{(X,Y)}(s,t)$. Finally, let $y'\neq y$ and $(s',t')= \nabla \Lambda^*_{(X,Y)}(x_0,y')$. Remark that $(x_0,y)= \nabla \Lambda_{(X,Y)}(s,t)$ and $(x_0,y')=\nabla \Lambda_{(X,Y)}(s',t')$ lead to $(s',t')\neq (s,t)$.  Therefore, by the strict convexity of  $\Lambda_{(X,Y)}$,
\begin{align*}
&\pscal{(x_0,y')-(x_0,y),\nabla \Lambda^*_{(X,Y)}(x_0,y')-\nabla \Lambda^*_{(X,Y)}(x_0,y)} \\
&\qquad \qquad =\pscal{\nabla \Lambda_{(X,Y)}(s',t')-\nabla \Lambda_{(X,Y)}(s,t),(s',t')-(s,t)}>0.
\end{align*}
Thus $\Lambda_{(X,Y)}^*(x_0,\cdot{})$ is strictly convex.  \qed
%Éventuellement à préciser :
%\begin{itemize}
%\item $m_N/N = 1/(1-\mu_N) \in R_{X^{(N)}} = \intervalleoo{0}{\infty}$ ;
%\item $1/(1-\mu) \in R_{X} = \intervalleoo{0}{\infty}$ ;
%\item $\ensavec{f'(u)}{u \leqslant 0} \supset \intervallefo{0}{\infty} \supset \intervallefo{0}{\Espe[Y]}$.
%\end{itemize}
%\end{proof}

%%%%%%%%%%%%%%%%%%%%%%%%%%%%%%%%%%
\subsection{Upper moderate deviations - Theorem \ref{th:umld_sparse} (i)}
%%%%%%%%%%%%%%%%%%%%%%%%%%%%%%%%%%

%\begin{proof}[Proof of \eqref{eq:umd_sparse} in Theorem \ref{th:umld_sparse}]
Analogously to \eqref{eq:prob}, one has 
\begin{align}
P_N \defeq \Prob(d_{m_N,n_N} - \Espe[d_{m_N,n_N}] \geqslant N^{\alpha} y)%\nonumber\\
%&= \Prob(\rva{T}{N} - \Espe[\rva{T}{N}\, | \, \rva{S}{N}=m_N] \geqslant N^{\alpha} y\, | \, \rva{S}{N}=m_N)\nonumber\\
&= \Prob(\rva{T}{N} - \Espe[\rva{T}{N}] \geqslant N^{\alpha} y_N\, | \, \rva{S}{N}=m_N)\label{eq:maj_conditionnement}\\
&\leqslant \frac{\Prob (\rva{T}{N} - \Espe[\rva{T}{N}] \geqslant N^{\alpha} y_N)}{\Prob (\rva{S}{N}=m_N)}\label{eq:maj_conditionnement2}
\end{align}
with $y_N \defeq  y + N^{-\alpha}( \Espe[T_N \, | \, S_N=m_N] - \Espe[T_N] ) \to y
$ by Proposition \ref{prop:moy_weak_array_cond_mob_cor_moderate}.
The upper bound then follows from Propositions \ref{prop:lognag} and \ref{prop:LLT_gd_dev}. 
As for the lower bound, using \eqref{eq:maj_conditionnement}, one has
\begin{align*}
P_N & \geqslant \Prob (\rva{T}{N} - \Espe[\rva{T}{N}] \geqslant N^{\alpha} y_N, \; \forall i, \rva[i]{Y}{N}<N^{\alpha/2}\, | \, \rva{S}{N}=m_N) \\
 & = \Prob (\rva{T}{N} - \Espe[\rva{T}{N}] \geqslant N^{\alpha} y_N \, | \, \forall i, \rva[i]{Y}{N}<N^{\alpha/2}, \; \rva{S}{N}=m_N)\Prob(\forall i, \rva[i]{Y}{N}<N^{\alpha/2}\, | \, \rva{S}{N}=m_N).
\end{align*}

On the one hand, one has 
\begin{equation} \label{eq:Sinf}
\Prob(\forall i, \rva[i]{Y}{N}<N^{\alpha/2}\, | \, \rva{S}{N}=m_N)
 = 1 - \Prob(\exists i, \rva[i]{Y}{N} \geqslant N^{\alpha/2}\, | \, \rva{S}{N}=m_N)
 \geqslant 1 - \frac{N \Prob(\rva[i]{Y}{N} \geqslant N^{\alpha/2})}{\Prob( \rva{S}{N}=m_N)} .
\end{equation}
Using Propositions \ref{prop:queue_Y} and \ref{prop:LLT_gd_dev}, we derive that $\Prob(\forall i, \rva[i]{Y}{N}<N^{\alpha/2}\, | \, \rva{S}{N}=m_N)\to 1$. On the other hand, let us turn to the minoration of $\Prob (\rva{T}{N} - \Espe[\rva{T}{N}] \geqslant N^{\alpha} y_N \, | \, \forall i, \rva[i]{Y}{N}<N^{\alpha/2}, \; \rva{S}{N}=m_N)$. 
In order to apply Gärtner-Ellis theorem, we follow the proof of \cite[Theorem 2.2]{TFC12} and we introduce
\begin{align*}
g_N(u)=\frac{1}{N^{2\alpha-1}}\log \Espe\bigl[e^{u (\rva{T}{N}-\Espe[\rva{T}{N}])/N^{1-\alpha}}\, \big| \, \forall i,\; \rva[i]{Y}{N}< N^{\alpha/2}, \; \rva{S}{N}=m_N\bigr].
\end{align*}
Write
\begin{align*}
\Espe&\bigl[e^{u (\rva{T}{N}-\Espe[\rva{T}{N}])/N^{1-\alpha}}\, \big| \, \forall i,\; \rva[i]{Y}{N}< N^{\alpha/2}, \; \rva{S}{N}=m_N\bigr]\\
&=\frac{\Espe\bigl[e^{u (\rva{T}{N}-\Espe[\rva{T}{N}])/N^{1-\alpha}}\indic_{\rva{S}{N}=m_N}\, \big| \, \forall i,\; \rva[i]{Y}{N}< N^{\alpha/2}\bigr]}{\Prob(\rva{S}{N}=m_N \, | \, \forall i, \rva[i]{Y}{N}<N^{\alpha/2})}\\
&= \frac{\Espe\bigl[e^{u (\rva{\Tinf}{N}-\Espe[\rva{T}{N}])/N^{1-\alpha}}\indic_{\rva{\Sinf}{N}=m_N} \bigr]}{\Prob(\rva{\Sinf}{N}=m_N)},
\end{align*}
where $\rva{\Sinf}{N}=\sum_{i=1}^{N} \rva[i]{\Xinf}{N}$, $\rva{\Tinf}{N}=\sum_{i=1}^{N} \rva[i]{\Yinf}{N}$, and  the random vectors $(\rva[i]{\Xinf}{N}, \rva[i]{\Yinf}{N})$ are independent, each distributed as $\mathcal L((\rva{X}{N},\rva{Y}{N})\, | \, \rva{Y}{N} < N^{\alpha/2})$. 
Then,
\begin{align*}
\Espe\bigl[e^{u (\rva{\Tinf}{N}-\Espe[\rva{T}{N}])/N^{1-\alpha}}\indic_{\rva{\Sinf}{N}=m_N} \bigr]
&=\frac{1}{2\pi}\int_{-\pi}^{\pi} e^{-ism_N}\Espe\bigl[e^{u (\rva{\Tinf}{N}-\Espe[\rva{T}{N}])/N^{1-\alpha}+is\rva{\Sinf}{N}} \bigr]ds\\
&=\frac{1}{2\pi}\int_{-\pi}^{\pi} e^{-ism_N}\Espe\bigl[e^{u (\rva{\Yinf}{N}-\Espe[\rva{Y}{N}])/N^{1-\alpha}+is\rva{\Xinf}{N}} \bigr]^{N}ds\\
&=e^{N\Lambda_{\rva{\Yinf}{N}-\Espe[\rva{Y}{N}]}(u/N^{1-\alpha})}\frac{1}{2\pi}\int_{-\pi}^{\pi} e^{-ism_N}\Espe\bigl[e^{is\rva{\hat X^u}{N}} \bigr]^{N}ds\\
&= e^{N\Lambda_{\rva{\Yinf}{N}-\Espe[\rva{Y}{N}]}(u/N^{1-\alpha})}\Prob(\rva{\hat S^u}{N}=m_N)
\end{align*}
where $\rva{\hat S^u}{N}$ stands for $\sum_{i=1}^{N} 
\rva[i]{\hat X}{N}^u$ and the random variables $\rva[i]{\hat X}{N}^u$ are independent copies of $\rva{\hat X}{N}^u$, the distribution of which is given by
\[
\Prob(\rva{\hat X}{N}^u=x)\defeq e^{-\Lambda_{\rva{\Yinf}{N}-\Espe[\rva{Y}{N}]} (u/N^{1-\alpha}) }\Espe\bigl[e^{u (\rva{\Yinf}{N}-\Espe[\rva{Y}{N}])/N^{1-\alpha}}\indic_{\rva{\Xinf}{N} =x}\bigr].
\]
Consequently,
\[
g_N(u) = N^{2-2\alpha} \Lambda_{\rva{\Yinf}{N}-\Espe[\rva{Y}{N}]}(u/N^{1-\alpha}) + \frac{1}{N^{2\alpha-1}} \log \Prob(\rva{\hat S^u}{N}=m_N) - \frac{1}{N^{2\alpha-1}} \log \Prob(\rva{\Sinf}{N}=m_N). 
\]

%\alert{Deux versions, la seconde colle plus à \cite{TFC12} :}
%
%\textbullet{} Since \alert{(calcul analogue fait en \eqref{eq:Sinf})}
%\[
%\Prob(\rva{\Sinf}{N}=m_N) = \frac{\Prob(\rva{S}{N} = m_N\ ;\ \forall i,\; \rva[i]{Y}{N}< N^{\alpha/2})}{\Prob(\forall i,\; \rva[i]{Y}{N}< N^{\alpha/2})} \geqslant \Prob(\rva{S}{N} = m_N) - N \Prob(\rva{Y}{N}\geqslant N^{\alpha/2})
%\]
%and applying Propositions \ref{prop:queue_Y} and \ref{prop:LLT_gd}, we get:
%\begin{align*}
%g_N(u)
% & = -N^{2-2\alpha}  \bigl( \Lambda_{\rva{\hat X}{N}^u}^*(m_N/N) - \Lambda_{\rva{\Yinf}{N}-\Espe[\rva{Y}{N}]}(u/N^{1-\alpha}) \bigr)(1 + o(1)) \\
% & = -N^{2-2\alpha} H_N(u/N^{1-\alpha}) (1 + o(1)) ,
%\end{align*}
%where
%\[
%H_N(t) = \sup \ensavec{s \frac{m_N}{N} - \Lambda_{(\rva{\Xinf}{N},\rva{\Yinf}{N}-\Espe[\rva{Y}{N}])}(s,t)}{s \in \R} .
%\]
%
%\textbullet{} 
So, using Lemma \ref{lem:LLT_gd_dev_bis} below, we get
\begin{align*}
g_N(u)
 & = - N^{2-2\alpha} \bigl( \Lambda_{\rva{\hat X}{N}^u}^*(m_N/N) - \Lambda_{\rva{\Xinf}{N}}^*(m_N/N) - \Lambda_{\rva{\Yinf}{N}-\Espe[\rva{Y}{N}]}(u/N^{1-\alpha}) \bigr) + O\Bigl( \frac{\log(N)}{N^{2 \alpha - 1}} \Bigr) \\
 & = -N^{2-2\alpha} \bigl(H_N(u/N^{1-\alpha}) - H_N(0)\bigr) + O\Bigl( \frac{\log(N)}{N^{2 \alpha - 1}} \Bigr) ,
\end{align*}
where
\[
H_N(t) = \sup \ensavec{s \frac{m_N}{N} - \Lambda_{(\rva{\Xinf}{N},\rva{\Yinf}{N}-\Espe[\rva{Y}{N}])}(s,t)}{s \in \R} .
\]

%\alert{À partir d'ici, tout reprendre. Notamment, $s(0) \neq 0$ et $H_N(0) \neq 0$, même si c'est le cas asymptotiquement.}

Applying the global version of the inverse function theorem to the function
$(s,t)\in \R\times\R\mapsto (t,s {m_N}/{N} - \Lambda_{(\rva{\Xinf}{N},\rva{\Yinf}{N}-\Espe[\rva{Y}{N}])}(s,t))$ and noting that $\partial_{s,s}\Lambda_{(\rva{\Xinf}{N},\rva{\Yinf}{N}-\Espe[\rva{Y}{N}])}(s,t)$ is nonzero since it is the variance of a non constant random variable, there is a unique maximizer $s_N(t)$ in the definition of $H_N(t)$. Moreover, the same algebraic computations as in \cite{TFC12} yield
\[
H_N'(0)=\Espe[\rva{Y}{N}]-\Espe[\rva{\Ytilde}{N}] \quad \text{and} \quad H_N''(0)=\frac{\Cov(\rva{\Ytilde}{N},\rva{\Xtilde}{N})^2}{\Var(\rva{\Xtilde}{N})}-\Var(\rva{\Ytilde}{N}),
\]
where the distribution of $(\rva{\Xtilde}{N},\rva{\Ytilde}{N})$ is given by
\begin{align}\label{def:XY_tilt}
\Prob(\rva{\Xtilde}{N}=x,\rva{\Ytilde}{N}=y)=e^{\tau_N x-\Lambda_{\rva{\Xinf}{N}}(\tau_N)}\Prob(\rva{\Xinf}{N}=x,\rva{\Yinf}{N}=y)
\end{align}
and $\tau_N$ is the unique solution of $\Lambda_{\rva{\Xinf}{N}}'(\tau_N)=m_N/N$. 
Then, 
%for $t_N = u/N^{1-\alpha} \to 0$,
\[
g_N(u) = -N^{1-\alpha} u H_N'(0)- \frac{u^2}{2} H_N''(0) -\frac{u^3}{6N^{1-\alpha}} H_N'''(z_N) + O\Bigl( \frac{\log(N)}{N^{2 \alpha - 1}} \Bigr)
\]
with $z_N \in \intervalleff{0}{u/N^{1-\alpha}}$. 
Using Remark \ref{rem:cv_espe_ytilde} below, one has $N^{1-\alpha}H_N'(0) \to 0$. By Lemma \ref{cor:cv_moments2} below, by \cite[Equation (4.31)]{Janson01a}, and by \eqref{def:sigma_mu}, we get
\[
H_N''(0) \to \frac{\Cov(X,Y)^2}{\Var(X)}-\Var(Y)=-\sigma^2(\mu).
\]
As in \cite{TFC12}, shortening $\Lambda_{(\rva{\Xinf}{N},\rva{\Yinf}{N}-\Espe[\rva{Y}{N}])}$ into $\Lambda$ and using obvious notation for partial derivatives, one has
\[
H_N'''(z_N) = \biggl( \biggl( \frac{\Lambda''_{s,t}}{\Lambda''_{s,s}} \biggr)^3 \Lambda'''_{s,s,s} - 3 \biggl( \frac{\Lambda''_{s,t}}{\Lambda''_{s,s}} \biggr)^2 \Lambda'''_{s,s,t} + 3 \frac{\Lambda''_{s,t}}{\Lambda''_{s,s}} \Lambda'''_{s,t,t} - \Lambda'''_{t,t,t} \biggr)(s_N(z_N),z_N).
\]
Let us prove that $s_N(z_N) \to 0$. The sequence of concave functions
\[
f_N(s)= s\frac{m_N}{N} - \Lambda_{(\rva{\Xinf}{N},\rva{\Yinf}{N}-\Espe[\rva{Y}{N}])}(s,z_N)
\]
converges pointwise to the strictly concave function $f(s)\defeq s(1-\mu)^{-1} -\Lambda_X(s)$. This fact follows from the uniform integrability of $\exp(s \rva{\Xinf}{N} + z_N (\rva{\Yinf}{N}-\Espe[\rva{Y}{N}]))$, which is a consequence of Lemma \ref{cor:cv_moments2} and the fact that $z_N\rva{\Yinf}{N}\leqslant u/N^{1-3\alpha/2}$ is bounded (remember that $\alpha < 2/3$).
Now the maximum of $f$ is attained at $0$. Let $\varepsilon\in \intervalleoo{0}{-\log(\lambda e)}$.  By the strict concavity of $f$, there exists $\eta>0$ such that 
$f(0)-\eta>\max(f(-\varepsilon),f(\varepsilon))+\eta$.
By \cite[Theorem 10.8]{Rockafellar70convex}, for all $N$ large enough, $\norme{f_N-f}_{\infty}<\eta$, where  
$\norme{\cdot{}}_{\infty}$ is the supremum norm over the compact set
 $\intervalleff{-\varepsilon}{\varepsilon}$. For those $N$, $s_N(z_N)\in \intervalleff{-\varepsilon}{\varepsilon}$.
Since $\varepsilon$ is arbitrary, we have proved that $s_N(z_N) \to 0$.

The uniform integrability of $
(\rva{\Xinf}{N})^p \abs{\rva{\Yinf}{N}-\Espe[\rva{Y}{N}]}^q \exp(s_N(z_N) \rva{\Xinf}{N} + z_N (\rva{\Yinf}{N}-\Espe[\rva{Y}{N}]))$
follows from the same arguments as before and 
the fact that $s_N(z_N) \to 0$, and Proposition \ref{prop:cv_moments} entails
\[
\Espe\bigl[ (\rva{\Xinf}{N})^p \abs{\rva{\Yinf}{N}-\Espe[\rva{Y}{N}]}^q e^{s_N(z_N) \rva{\Xinf}{N} + z_N (\rva{\Yinf}{N}-\Espe[\rva{Y}{N}]) }\bigr] \to \Espe\bigl[X^p \abs{Y-\Espe[Y]}^q\bigr] .
\]
Therefore, $H_N'''(z_N)$ is bounded, whence
$g_N(u) \to u^2 \sigma^2(\mu)/2$ and \eqref{eq:umd_sparse} follows. \qed
%\end{proof}

\begin{lem}\label{lem:LLT_gd_dev_bis}
Let $\rva{\check X}{N}$ be $\rva{\Xinf}{N}$ or $\rva{\hat X}{N}^u$. Denoting by $\rva{\check S}{N}\defeq\rva[1]{\check X}{N}+\dots \rva[N]{\check X}{N}$, 
 we have, for any sequence of integers $(k_N)_{N \geqslant 1}$ such that $\lim k_N/N\in \intervalleoo{1}{\infty}$,
\[
\log \Prob(\rva{\check S}{N} = k_N) = - N \Lambda_{\rva{\check X}{N}}^*(k_N/N) + O(\log(N)) .
\]
\end{lem}

\begin{proof}[Proof of Lemma \ref{lem:LLT_gd_dev_bis}]
We just check that we can apply \cite[Lemma 3.3]{TFC12} to the sequences $(\rva{\Xinf}{N})_{N \geqslant 1}$ and $(\rva{\hat X}{N}^u)_{N \geqslant 1}$. The conclusion follows since, in this case, $m = 1$, $b = 0$ and $c_{n, m, b} = 1$.

\textbullet{} First, $\Ima(\Lambda_{\rva{\Xinf}{N}}')=\intervalleoo{1}{\infty}=\Ima(\Lambda_X')$ so, for all $N$ large enough, $k_N/N\in \Ima(\Lambda_{\rva{\Xinf}{N}}')$ and $\lim k_N/N\in \Ima(\Lambda_X')$. Secondly,
\[
\Espe[e^{s \rva{\Xinf}{N}}]
 = \sum_{x \geqslant 1} e^{s x} \frac{\Prob(\rva{X}{N} = x,\ \rva{Y}{N} < N^{\alpha / 2})}{\Prob(\rva{Y}{N} < N^{\alpha / 2})}
 \leqslant \sum_{x \geqslant 1} e^{s x} \frac{\Prob(\rva{X}{N} = x)}{\Prob(\rva{Y}{N} < N^{\alpha / 2})}
 = \frac{\Espe[e^{s \rva{X}{N}}]}{\Prob(\rva{Y}{N} < N^{\alpha / 2})} ,
\]
so
\[
\interieur(\dom(\Lambda_{\rva{\Xinf}{N}}))
 \supset \interieur(\dom(\Lambda_{\rva{X}{N}}))
 = \intervalleoo{-\infty}{-\log(\lambda_N e)} .
\]
Since $\interieur(\dom(\Lambda_X))=\intervalleoo{-\infty}{-\log(\lambda e)}$ and $\lambda_N \to \lambda$, assumption 1.\ holds for $N$ large enough. Finally, it remains to check that assumption 2.\ is satisfied. For $s<-\log(\lambda e)$ and $t\in\R$, we have
\begin{align*}
&\abs{\Espe\bigl[e^{(s+it)\rva{\Xinf}{N}}\bigr]-\Espe\bigl[e^{(s+it)X}\bigr]}\\
&\leqslant \frac{\abs{\Espe\bigl[e^{(s+it)\rva{X}{N}}\bigr]-\Espe\bigl[e^{(s+it)X}\bigr]}+\Prob(\rva{Y}{N}> N^{1/2})\Espe\bigl[e^{sX}\bigr]+\Espe\bigl[e^{s\rva{X}{N}}\indic_{\rva{Y}{N}> N^{1/2}}\bigr]}{\Prob(\rva{Y}{N}\leqslant  N^{1/2})}.
\end{align*}
%Now,
%\[
%\abs{\Espe[e^{i s \rva{X}{N}}] - \Espe[e^{i s X}]}
% \leqslant \sum_{x=1}^\infty \abs{\Prob(\rva{X}{N} = x) - \Prob(X = x)}
% \to 0
%\]
Now, for $s'\in \intervalleoo{s}{-\log(\lambda e)}$, using Hölder's inequality in the third line below,
\begin{align*}
&\sup_{t\in \R}\abs{\Espe[e^{(s+it) \rva{X}{N}}] - \Espe[e^{ (s+it) X}]}\\
& \hspace{2cm}
\leqslant \sum_{x=1}^\infty e^{sx} \abs{\Prob(\rva{X}{N} = x) - \Prob(X = x)}\\
 & \hspace{2cm} = \sum_{x=1}^\infty e^{sx} \abs{\Prob(\rva{X}{N} = x) - \Prob(X = x)}^{s/s'}\cdot{} \abs{\Prob(\rva{X}{N} = x) - \Prob(X = x)}^{1-s/s'}\\
 & \hspace{2cm} \leqslant \Bigl(\sum_{x=1}^\infty e^{s'x} \abs{\Prob(\rva{X}{N} = x) - \Prob(X = x)}\Bigr)^{s/s'}\Bigl(\sum_{x=1}^\infty \abs{\Prob(\rva{X}{N} = x) - \Prob(X = x)}\Bigr)^{1-s/s'}\\
 & \hspace{2cm} \leqslant \Bigl(\Espe[e^{s'\rva{X}{N}}]+\Espe[e^{s'X}]\Bigr)^{s/s'}\Bigl(\sum_{x=1}^\infty \abs{\Prob(\rva{X}{N} = x) - \Prob(X = x)}\Bigr)^{1-s/s'}\\
 & \hspace{2cm} \to 0 ,
\end{align*}
by Proposition \ref{prop:cv_moments} (for discrete random variables, the convergence in distribution is equivalent to the convergence in total variation); hence, the first term of the numerator converges to $0$ uniformly in $t \in \R$. So does also the second one by Proposition \ref{prop:queue_Y}. Finally, by the same arguments, for $s'\in \intervalleoo{s}{-\log(\lambda e)}$, 
\begin{align*}
\Espe\bigl[e^{s\rva{X}{N}}\indic_{\rva{Y}{N}> N^{1/2}}\bigr]
 & \leqslant \Espe\bigl[e^{s' \rva{X}{N}}\bigr]^{s/s'} \Prob(\rva{Y}{N}> N^{1/2})^{1-s/s'}\to 0,
\end{align*}
leading to the required result.
%Finally, by the same arguments and Hölder's inequality, with $1/p+1/q=1$ and $ps<-\log(\lambda e)$, 
%\begin{align*}
%\Espe\bigl[e^{s\rva{X}{N}}\indic_{\rva{Y}{N}> N^{1/2}}\bigr]& \leqslant
%\Espe\bigl[e^{p s\rva{X}{N}}\bigr]^{1/p}\Prob(\rva{Y}{N}> N^{1/2})^{1/q}\to 0,
%\end{align*}
%leading to the required result.

\textbullet{} As before, $\Ima(\Lambda_{\rva{\hat X}{N}^u}')=\intervalleoo{1}{\infty}=\Ima(\Lambda_X')$ so, for all $N$ large enough, $k_N/N\in \Ima(\Lambda_{\rva{\hat X}{N}^u}')$ and $\lim k_N/N\in \Ima(\Lambda_X')$. Since $\rva{\Yinf}{N}$ is bounded, $\interieur(\dom(\Lambda_{\rva{\hat X}{N}^u})) \supset \interieur(\dom(\Lambda_{\rva{\Xinf}{N}}))$, so assumption 1.\ holds. Finally, using the definition of $\rva{\hat X}{N}^u$, one gets
\begin{align*}
\Espe[e^{(s+it)\rva{\hat X}{N}^u}]
=e^{-\Lambda_{\rva{\Yinf}{N}-\Espe[\rva{Y}{N}]}(u/N^{1-\alpha})}
\Espe\left[e^{(s+it)\rva{\Xinf}{N}+\frac{u}{N^{1-\alpha}} (\rva{\Yinf}{N}-\Espe[\rva{Y}{N}])}\right]
\end{align*}
that converges to $\Espe[e^{(s+it)X}]$ uniformly in $t\in \R$ by similar arguments, and assumption 2.\ is satisfied.
\end{proof}

\begin{lem} \label{lem:cv_tau}
Let $\tau_N$ be the unique solution of $\Lambda_{\rva{\Xinf}{N}}'(\tau_N)=m_N/N$. There exists $c > 0$ such that, for all $N$ large enough, $\abs{\tau_N} \leqslant e^{-cN^{\alpha/4}}$.
\end{lem}

\begin{proof}[Proof of Lemma \ref{lem:cv_tau}]
First, for all $s < -\log(\lambda e)$,
\begin{equation} \label{eq:ecart_lambda_prime}
\abs{\Lambda_{\rva{\Xinf}{N}}'(s) - \Lambda_{\rva{X}{N}}'(s)}
= \abs{\frac{\Espe\bigl[\rva{X}{N}e^{s \rva{X}{N}}\bigr] \Espe\bigl[e^{s \rva{X}{N}} \indic_{\rva{Y}{N} \geqslant N^{\alpha/2}}\bigr]}{\Espe\bigl[e^{s \rva{X}{N}}\bigr] \Espe\bigl[e^{s \rva{X}{N}} \indic_{\rva{Y}{N} < N^{\alpha/2}}\bigr]} - \frac{\Espe\bigl[\rva{X}{N} e^{s \rva{X}{N}} \indic_{\rva{Y}{N} \geqslant N^{\alpha/2}}\bigr]}{\Espe\bigl[e^{s \rva{X}{N}} \indic_{\rva{Y}{N} < N^{\alpha/2}}\bigr]}} \leqslant e^{-c_1 N^{\alpha/4}}
\end{equation}
for some constant $c_1 > 0$ (independent of $s$ and $N$), using Hölder's inequality and Propositions \ref{prop:queue_Y} and \ref{prop:cv_moments}. Now, write
\[
\Lambda_{\rva{X}{N}}'(s) = \Espe[\rva{X}{N}] + s \Var(\rva{X}{N}) + \frac{s^2}{2} \Lambda_{\rva{X}{N}}'''(t)
\]
with $t$ between $0$ and $s$. Using Proposition \ref{prop:cv_moments}, there exists $s_0 > 0$ such that, for all $s \in \intervalleff{-s_0}{s_0}$ and for all $N$ large enough,
\begin{equation} \label{eq:dl_lambda_prime}
\abs{\Lambda_{\rva{X}{N}}'(s) - \Espe[\rva{X}{N}] - s \Var(\rva{X}{N})} \leqslant \frac{\abs{s} \Var(\rva{X}{N})}{2}.
\end{equation}
Since $\Lambda_{\rva{\Xinf}{N}}'(\tau_N) = \Espe[\rva{X}{N}]$,  \eqref{eq:ecart_lambda_prime} and \eqref{eq:dl_lambda_prime} yield
$\abs{\tau_N} \leqslant 2 e^{-c_1 N^{\alpha/4}}/\Var(\rva{X}{N})$, hence the desired result since $\Var(\rva{X}{N})\to \Var(X)>0$.
\end{proof}

\begin{lem}\label{cor:cv_moments2}
$(\rva{\Xinf}{N},\rva{\Yinf}{N})_{N \geqslant 1}$ and $(\rva{\Xtilde}{N},\rva{\Ytilde}{N})_{N \geqslant 1}$ converge to  $(X,Y)$ in distribution and with all mixed moments of the type $\Espe[\rva{\check{X}}{N}^p \rva{\check{Y}}{N}^qe^{s\rva{\check{X}}{N}}]$, where $p\geqslant 0$, $q\geqslant 0$, $s<-\log(\lambda e)$, and $\check{X}$ (resp.\ $\check{Y}$) stands for $\Xinf$ or $\Xtilde$ (resp.\ $\Yinf$ or $\Ytilde$).
\end{lem}

\begin{proof}[Proof of Lemma \ref{cor:cv_moments2}]
Following the proof of Proposition \ref{prop:cv_moments}, we prove separately that  $(\rva{\Xinf}{N})_{N \geqslant 1}$ and $(\rva{\Xtilde}{N})_{N \geqslant 1}$ converge to  $X$ in distribution and with all moments and $(\rva{\Yinf}{N})_{N \geqslant 1}$ and $(\rva{\Ytilde}{N})_{N \geqslant 1}$ converge to  $Y$ in distribution and with all moments. 

Let us prove the convergence of $(\rva{\Xinf}{N})_{N \geqslant 1}$. Let $f\colon \R\mapsto\R$ be a bounded measurable function. By Proposition  \ref{prop:cv_moments}, it suffices to prove that $\Espe[f(\rva{\Xinf}{N})]-\Espe[f(\rva{X}{N})]|$ converges to $0$. One has
\begin{align*}
\abs{\Espe[f(\rva{\Xinf}{N})]-\Espe[f(\rva{X}{N})]}
 & = \abs{\Espe[f(\rva{X}{N}) (\Prob(\rva{Y}{N} \geqslant N^{\alpha/2}) - \indic_{\rva{Y}{N} \geqslant N^{\alpha/2}})]} \Prob(\rva{Y}{N} < N^{\alpha/2})^{-1} \\
 & \leqslant 2 \norme{f}_\infty \Prob(\rva{Y}{N} \geqslant N^{\alpha/2}) \Prob(\rva{Y}{N} < N^{\alpha/2})^{-1} .
\end{align*}
%\begin{align*}
%|\Espe[f(\rva{\Xinf}{N})]-\Espe[f(\rva{X}{N})]|
%&\leqslant (\abs{\Espe[f(\rva{X}{N})]}+\Espe[f(\rva{X}{N})^2]^{1/2})\frac{\Prob(\rva{Y}{N} \geqslant N^{\alpha/2})}{\Prob(\rva{Y}{N} < N^{\alpha/2})}.
%\end{align*}
The result follows from Proposition \ref{prop:queue_Y}. To prove the convergence of the moments of $(\rva{\Xinf}{N})_{N \geqslant 1}$ to  those of $X$, it suffices to show that, for all $r > 0$, $(\Espe[(\rva{\Xinf}{N})^r])_{N \geqslant 1}$ is bounded (see, e.g., \cite[Example 2.21]{van1998asymptotic}) and to conclude by uniform integrability. We have
\begin{align*}
\Espe[(\rva{\Xinf}{N})^r]
 &   = \Espe[\rva{X}{N}^{r}\, | \, \rva{Y}{N}<N^{1/2}] \leqslant \frac{\Espe[\rva{X}{N}^{r}]}{\Prob (\rva{Y}{N}<N^{1/2})},
\end{align*}
so $\limsup \Espe[(\rva{\Xinf}{N})^r]$ is finite by Propositions \ref{prop:queue_Y} and \ref{prop:cv_moments}. The same calculation leads to the convergence of $(\Espe[\exp(s\rva{\Xinf}{N})])_{N\geqslant 1}$ and the same lines yield the convergence in distribution of $(\rva{\Yinf}{N})_{N\geqslant 1}$ to  $Y$ and with all moments.

Let us consider $(\rva{\Xtilde}{N})_{N\geqslant 1}$. Let $f\colon \R\mapsto\R$ be a bounded measurable function. One has
\begin{align}
|\Espe[f(\rva{\Xtilde}{N})]-\Espe[f(\rva{\Xinf}{N})]|
&\leqslant 
\frac{\left|\Espe\left[f(\rva{\Xinf}{N})\bigl(e^{\tau_N \rva{\Xinf}{N}}-\Espe\bigl[e^{\tau_N \rva{\Xinf}{N}}\bigr]\bigr)\right]\right|}{\Espe\left[e^{\tau_N \rva{\Xinf}{N}}\right]}
\leqslant \norme{f}_{\infty}
\frac{\Espe\left[\left|e^{\tau_N \rva{\Xinf}{N}}-\Espe\bigl[e^{\tau_N \rva{\Xinf}{N}}\bigr]\right|\right]}{\Espe\left[e^{\tau_N \rva{\Xinf}{N}}\right]}.\label{eq:decomp}
\end{align}
Using the convergence of $(\Espe[\exp(s\rva{\Xinf}{N})])_{N\geqslant 1}$ together with Lemma \ref{lem:cv_tau}, Slutsky's Lemma, and uniform integrability, we get $\Espe\left[e^{\tau_N \rva{\Xinf}{N}}\right]\to 1$. Similarly, since 
\[
\left|e^{\tau_N \rva{\Xinf}{N}}-\Espe\bigl[e^{\tau_N \rva{\Xinf}{N}}\bigr]\right|^2=e^{2\tau_N\rva{\Xinf}{N}}-2 e^{\tau_N\rva{\Xinf}{N}}\Espe\bigl[e^{\tau_N \rva{\Xinf}{N}}\bigr] + \Espe\bigl[e^{\tau_N \rva{\Xinf}{N}}\bigr]^2,
\]
the numerator in the right-hand side of \eqref{eq:decomp} converges to $0$.
Hence $(\rva{\Xtilde}{N})_{N\geqslant 1}$ converges in distribution to $X$ and similar arguments show the convergence of all moments of $(\rva{\Xtilde}{N})_{N\geqslant 1}$.  
The same lines lead to the convergence in distribution of $(\rva{\Ytilde}{N})_{N\geqslant 1}$ to  $Y$ and with all moments.
\end{proof}

\begin{rem}\label{rem:cv_espe_ytilde}
Following the same lines as in the proof of Lemma \ref{cor:cv_moments2} for $f=\id_{\R}$ and using Proposition \ref{prop:cv_moments} instead of the fact that $f$ is bounded and the fact that $\tau_N$ converges exponentially rapidly to $0$ (by Lemma \ref{lem:cv_tau}), we get that $\Espe[\rva{\Ytilde}{N}]-\Espe[\rva{Y}{N}]$ converges exponentially rapidly to $0$. 
\end{rem}

We notice that the rate function in the upper moderate deviations $y^2/(2\sigma^2(\mu))$ in \eqref{eq:umd_sparse} depends on the conditioning on $\{\rva{S}{N}=m_N\}$. As opposed to this, the conditioning does not influence the expression of the rate function in the upper large deviations, due to the fact that the random variables $\rva{Y}{N}$ are heavy-tailed. As a consequence, our proof of \eqref{eq:uld_sparse} in Section \ref{ssec:proof_umld_sparse} does not mimic that of \cite[Theorem 2.1]{TFC12} but is rather inspired by that of \cite[Theorem 5]{Nagaev69-1}.

%%%%%%%%%%%%%%%%%%%%%%%%%%%%%%%%%%
\subsection{Upper intermediate deviations - Theorem \ref{th:umld_sparse} (ii)}
%%%%%%%%%%%%%%%%%%%%%%%%%%%%%%%%%%

%Remind that
%\begin{align*}
%P_N=\Prob(d_{m_N,n_N} - \Espe[d_{m_N,n_N}] \geqslant N^{\alpha} y)
%&= \Prob(\rva{T}{N} - \Espe[\rva{T}{N}] \geqslant N^{\alpha} y_N\, | \, \rva{S}{N}=m_N).
%\end{align*}

%\begin{proof}[Proof of \eqref{eq:uid_sparse} in Theorem \ref{th:umld_sparse}]
Here $\alpha = 2/3$. The upper bound comes from \eqref{eq:maj_conditionnement2} and Propositions \ref{prop:LLT_gd_dev} and \ref{prop:lognag} (ii).
Let us turn to the lower bound. 
We assume that the infimum in the right-hand side of \eqref{eq:queue_Y} is attained at $\delta_0$.
Let $z>0$ and 
$l_N\defeq \lceil (N^\alpha z/\delta_0)^{1/2}\rceil$. By \eqref{eq:maj_conditionnement}, we have
\begin{align*}
P_N & \geqslant \Prob( \rva{T}{N} - \Espe[\rva{T}{N}] \geqslant N^\alpha y_N,  \;  \rva{S}{N}= m_N ) \\
 & \geqslant \Prob\big( T_N - \Espe[T_N] \geqslant N^\alpha y_N,  \;  S_N= m_N,  \;  \rva[N]{Y}{N} - \Espe[\rva[N]{Y}{N}] \geqslant N^\alpha z \big) \\
&\geqslant  \Prob(\rva[N-1]{T}{N} - \Espe[\rva[N-1]{T}{N}]\geqslant N^\alpha (y_N-z)\vert \rva[N-1]{S}{N}= m_N-l_N )\Prob(\rva[N-1]{S}{N}= m_N-l_N )\\
& \qquad \qquad \qquad \Prob(\rva[N]{Y}{N}-\Espe[\rva[N]{Y}{N}] \geqslant N^\alpha z,  \;  \rva[N]{X}{N}=l_N ).
\end{align*}
By Proposition \ref{prop:queue_couple}, 
\begin{align*} %\label{P3_inter}
\liminf \frac{1}{N^{\alpha/2}} \log \Prob(\rva{Y}{N}-\Espe[\rva{Y}{N}] \geqslant N^\alpha z,  \;  \rva{X}{N}=l_N ) & = - q(\mu)z^{1/2}
\end{align*}
and, by Proposition \ref{prop:LLT_gd_dev}, we derive that
\begin{align*}
\frac{1}{N^{\alpha/2}}\log  \Prob( \rva[N-1]{S}{N} = m_N-l_N )
 & = -N^{1-\alpha/2} \Lambda_{\rva{X}{N}}^* \Bigl(\frac{m_N-l_N}{N}\Bigr) + O\Bigl( \frac{\log(N)}{N^{1 - \alpha / 2}} \Bigr).
\end{align*}
Let us prove that
 \begin{align}\label{eq:cvLambda*N}
 \Lambda_{\rva{X}{N}}^*\Bigl(\frac{m_N-l_N}{N}\Bigr) = \frac{1}{2 \sigma_X^2} \Bigl(\frac {l_N}{N}\Bigr)^2 + o\Bigl(\Bigl(\frac {l_N}{N}\Bigr)^2\Bigr).
 \end{align}
We have
\begin{align*}
\Lambda_{\rva{X}{N}}^* \Bigl(\frac{m_N-l_N}{N}\Bigr) &= \frac 12 \Bigl(\frac {l_N}{N}\Bigr)^2 \frac{1}{\ssdd{\rva{X}{N}}^2}-\frac 16 \Bigl(\frac{l_N}{N}\Bigr)^3(\Lambda_{\rva{X}{N}}^*)''' (c_N)
\end{align*}
where $c_N\in \intervalleff{\frac{m_N-l_N}{N}}{\frac{m_N}{N}}$ converges to $(1-\mu)^{-1}$. Then direct computations give, for all $c$,
\begin{align*}
(\Lambda_{\rva{X}{N}}^*)''' (c) = \frac{\Lambda_{\rva{X}{N}}''' (s_N(c))}{\Lambda_{\rva{X}{N}}'' (s_N(c))} \quad \text{and} \quad 
(\Lambda_X^*)''' (c) = \frac{\Lambda_X''' (s(c))}{\Lambda_X'' (s(c))} 
\end{align*}
where $s_N(c)$ (resp.\ $s(c)$) is the unique solution of $\Lambda_{\rva{X}{N}}'(s_N(c))=c$ (resp.\ $\Lambda_X'(s(c))=c$). Let us prove that $s_N(c)\to s(c)$. Since $\Lambda_X'' (s(c))=2\delta>0$, there exists $\alpha>0$ such that 
 $\Lambda_X'' >\delta$ over $V=\intervalleff{s(c)-\alpha}{s(c)+\alpha}$. Let $\varepsilon\in \intervalleoo{0}{\alpha}$. For $N$ large enough and $s \in V$, $\abs{\Lambda_{\rva{X}{N}}'(s)-\Lambda_X'(s)}\leqslant \delta\varepsilon$, by Proposition \ref{prop:cv_moments} and the uniform convergence of power series on compact subsets of the domain of convergence. Then 
\begin{align*}
\Lambda_{\rva{X}{N}}'(s(c)-\varepsilon)\leqslant \Lambda_X'(s(c)-\varepsilon)+\delta \varepsilon \leqslant \Lambda_X'(s(c))=c \leqslant \Lambda_X'(s(c)+\varepsilon)-\delta \varepsilon \leqslant  \Lambda_{\rva{X}{N}}'(s(c)+\varepsilon).
\end{align*} 
Since $\Lambda_{\rva{X}{N}}'$ is increasing, we deduce that $s_N(c)\in \intervalleff{s(c)-\varepsilon}{s(c)+\varepsilon}$. So $s_N(c)\to s(c)$ as announced. 
 Using Proposition \ref{prop:cv_moments} and the uniform convergence of power series on compact subsets of the domain of convergence, we conclude that
$
(\Lambda_{\rva{X}{N}}^*)''' (c_N) \to 
(\Lambda_X^*)''' (1/(1-\mu))$.
 Using Proposition \ref{prop:cv_moments} again, we get \eqref{eq:cvLambda*N}. Finally,
\begin{align*}
\frac{1}{N^{\alpha/2}}\log  \Prob( \rva[N-1]{S}{N} = m_N-l_N )
  & = -\frac{N^{1-\alpha/2} }{2\ssdd{X}^2}\Bigl(\frac{l_N}{N}\Bigr)^2 + O\Bigl( \frac{\log(N)}{N^{1 - \alpha / 2}} \Bigr),
\end{align*}
which converges to $0$ since $\alpha<2$. %and by Proposition \ref{prop:LLT_gd_dev}, 
%\begin{align*} %\label{P2_inter}
%\liminf \frac{1}{N^{\alpha/2}} \log \Prob(\rva[N-1]{S}{N}= m_N-l_N )  \to 0.
%\end{align*}
As for the minoration of the remaining term, we follow the same lines as in the proof of the lower bound in Theorem \ref{th:umld_sparse} (i) which remains valid for $\alpha=2/3$ and $\Espe[\rva[N-1]{S}{N}]=m_N+O(N^{1/2})$. Hence,
\begin{align*} %\label{P1_inter}
\liminf \frac{1}{N^{\alpha/2}} \log \Prob(\rva[N-1]{T}{N} - \Espe[\rva[N-1]{T}{N}]\geqslant N^\alpha (y_N-z)\, \vert \, \rva[N-1]{S}{N}= m_N-l_N ) = -\frac{(y-z)^2}{2\sigma^2(\mu)}.
\end{align*}
Optimizing in $z = (1 - t) y$ with $t\in \intervalleoo{0}{1}$ leads to \eqref{eq:uid_sparse}. \qed
%\end{proof}

%%%%%%%%%%%%%%%%%%%%%%%%%%%%%%%%%%
\subsection{Upper large deviations - Theorem \ref{th:umld_sparse} (iii)}\label{ssec:proof_umld_sparse}
%%%%%%%%%%%%%%%%%%%%%%%%%%%%%%%%%%

%We study the asymptotic behavior of the quantity $\Prob(d_{m_N,n_N} - \Espe[d_{m_N,n_N}] \geqslant N^{\alpha} y)$ for $\alpha>2/3$. Remind that the total displacement $d_{m_N,n_N}$ is distributed as the conditional distribution of $T_N$ given $S_N=m_N$. Notice that $\Espe[S_N] = N \Espe[\rva{X}{N}] = m_N$. Now let
%\begin{align*}
%P_N & \defeq \Prob(d_{m_N,n_N} - \Espe[d_{m_N,n_N}] \geqslant N^\alpha y) \\
% &= \Prob(T_N - \Espe[T_N| S_N=m_N] \geqslant N^\alpha y \, | \, S_N=m_N) \\
% &= \Prob(T_N - \Espe[T_N] \geqslant N^\alpha y_N, \;  S_N=m_N)/\Prob(S_N=m_N)
%\end{align*}
%where $y_N \defeq y + N^{-\alpha} (\Espe[T_N\, | \, S_N=m_N] - \Espe[T_N])$. Proposition \ref{prop:moy_weak_array_cond_mob_cor_moderate} entails $y_N \to y$. 
%
%\medskip

%\begin{proof}[Proof of \eqref{eq:uld_sparse} in Theorem \ref{th:umld_sparse}]
Using \eqref{eq:maj_conditionnement2}, the upper bound
\[
\limsup_{N \to \infty} \frac{1}{N^{\alpha/2}} \log(P_N) \leqslant \limsup_{N \to \infty} \frac{1}{N^{\alpha/2}} \log \biggl( \frac{\Prob(T_N - \Espe[T_N] \geqslant N^\alpha y_N)}{\Prob(S_N=m_N)} \biggr) = -q(\mu) y^{1/2},
\]
follows from Propositions \ref{prop:queue_Y}, \ref{prop:LLT_gd_dev}, and \ref{prop:lognag} (iii). For the lower bound, assume that the infimum in the right-hand side of \eqref{eq:queue_Y} is attained at $\delta_0$. Let $\varepsilon>0$ and $l_N\defeq \lceil (N^\alpha(y_N + \varepsilon)/\delta_0)^{1/2}\rceil$. By \eqref{eq:maj_conditionnement}, we have
\begin{align*}
P_N & \geqslant \Prob( T_N - \Espe[T_N] \geqslant N^\alpha y_N,  \;  S_N= m_N ) \\
 & \geqslant \Prob\big( T_N - \Espe[T_N] \geqslant N^\alpha y_N,  \;  S_N= m_N,  \;  \rva[N]{Y}{N} - \Espe[\rva[N]{Y}{N}] \geqslant N^\alpha(y_N+\varepsilon) \big) \\
 & \geqslant \Prob\big( \rva[N-1]{T}{N} - \Espe[\rva[N-1]{T}{N}] \geqslant - N^\alpha\varepsilon,  \; \rva[N-1]{S}{N} = m_N-l_N\big) \\
 & \qquad \qquad \qquad \Prob\big(\rva{Y}{N} - \Espe[\rva{Y}{N}] \geqslant N^\alpha(y_N+\varepsilon),  \;  \rva{X}{N}=l_N \big) \\
 & \eqdef P_{N,1} P_{N,2}.
\end{align*}
Applying Proposition \ref{prop:queue_couple}, one gets
\begin{align} \label{P2}
\lim_{N \to \infty} \frac{1}{N^{\alpha/2}} \log(P_{N,2}) 
%& = \liminf_{N \to \infty} \frac{1}{N^{\alpha/2}} \log ( \Prob( d_{l_N, l_N-1} \geqslant N^\alpha(y_N + \varepsilon) + \Espe[\rva{Y}{N}] ) \Prob( \rva{X}{N} = l_N )) \nonumber\\
% & = \liminf_{N \to \infty} \frac{1}{N^{\alpha/2}} \log ( \Prob( d_{l_N, l_N} \geqslant N^\alpha(y_N + \varepsilon) + \Espe[\rva{Y}{N}] ) \Prob( \rva{X}{N} = l_N ))\nonumber\\
% & \geqslant \liminf_{N \to \infty} \sqrt{\frac{y_N}{\delta_0}} \left(\frac{1}{l_N}\log ( \Prob( d_{l_N, l_N} \geqslant \delta_0 l_N^2) \Prob( \rva{X}{N} = l_N))\right)\nonumber\\
& = -\sqrt{\frac{y+\varepsilon}{\delta_0}} \left(\kappa(\mu) + J(\delta_0)\right) \underset{\varepsilon\to 0}{\to} - q(\mu)y^{1/2}.
\end{align}
Let us turn to the minoration of $P_{N,1}$. One has
\begin{align*}
P_{N,1}
% & = \Prob( \rva[N-1]{T}{N}- \Espe[\rva[N-1]{T}{N}] \geqslant -N^{\alpha}\varepsilon,  \;  \rva[N-1]{S}{N} = m_N-l_N) \\
 & \geqslant \Prob( \rva[N-1]{S}{N} = m_N-l_N ) - \Prob( \rva[N-1]{T}{N}- \Espe[\rva[N-1]{T}{N}] < -N^\alpha\varepsilon ).
\end{align*}
Applying Proposition \ref{prop:LLT_gd_dev} and \eqref{eq:cvLambda*N}, we derive that
\begin{align*}
\frac{1}{N^{\alpha/2}}\log  \Prob( \rva[N-1]{S}{N} = m_N-l_N )
 & = -N^{1-\alpha/2} \Lambda_{\rva{X}{N}}^* \Bigl(\frac{m_N-l_N}{N}\Bigr) + O\Bigl( \frac{\log(N)}{N^{1 - \alpha / 2}} \Bigr) \\
 & = -\frac{N^{1-\alpha/2} }{2\ssdd{X}^2}\Bigl(\frac{l_N}{N}\Bigr)^2 + O\Bigl( \frac{\log(N)}{N^{1 - \alpha / 2}} \Bigr),
% & = O\Bigl( \frac{\log(N)}{N^{1 - \alpha / 2}} \Bigr) ,
\end{align*}
%\[
%\frac{1}{N^{\alpha/2}}\log  \Prob( \rva[N-1]{S}{N} = m_N-l_N )
%\sim  -N^{1-\alpha/2} \Lambda_{\rva{X}{N}}^* \Bigl(\frac{m_N-l_N}{N}\Bigr) 
%\sim -\frac{N^{1-\alpha/2} }{2\ssdd{X}^2}\Bigl(\frac{l_N}{N}\Bigr)^2 =O \Bigl(\frac{1}{N^{1-\alpha/2}}\Bigr)
%\]
which converges to $0$ since $\alpha<2$. Applying a unilateral version of  \cite[Theorem 2.2]{TFC12}, we get
\[
\frac{1}{N^{\alpha/2}}\log  \Prob( \rva[N-1]{T}{N}- \Espe[\rva[N-1]{T}{N}] < -N^\alpha\varepsilon)  \begin{cases} 
\sim -c_{\varepsilon}N^{3\alpha/2-1}\to -\infty & \text{if $\alpha\in \intervalleof{2/3}{1}$}\\
=-\infty & \text{if $\alpha>1$}
\end{cases}
\]
for some $c_{\varepsilon}>0$. Eventually,
$N^{-\alpha/2} \log(P_{N,1}) \to 0$,
which leads, together with \eqref{P2}, to
\[
\pushQED{\qed} 
\liminf \frac{1}{N^{\alpha/2}} \log(P_N) \geqslant - q(\mu)y^{1/2}.\qedhere
\popQED
\]
%\end{proof}

%%%%%%%%%%%%%%%%%%%%%%%%%%%%%%%%%%
\subsection{Upper large deviations for \texorpdfstring{$\alpha=2$}{f} - Theorem \ref{th:umld_sparse} (iv)}
%%%%%%%%%%%%%%%%%%%%%%%%%%%%%%%%%%

%Let $y > 0$ and let $m_N \sim Nx$ with $x > 1$. Actually, we only need the case $x = (1-\mu)^{-1}$. Assume that $\mu$ satisfies:
%\begin{equation} \label{eq:mu}
%x-\sqrt{2y} \leqslant (1-\mu)^{-1} = \Espe[X],
%\end{equation}
%where $X$ is the limit of $(\rva{X}{N})_N$ defined in Proposition \ref{prop:cv_moments} (that obviously  holds when $x= (1-\mu)^{-1}$). 
%
%Then, using \eqref{eq:maj_conditionnement} and  Proposition \ref{prop:LLT_gd_dev} and reminding $\Espe[\rva{S}{N}] = N \Espe[\rva{X}{N}] = m_N$, it suffices to prove
%\begin{align}\label{eq:ud_sparse_particulier}
%\frac{1}{N} \log \Prob(\rva{S}{N} = m_N,  \; \rva{T}{N} \geqslant N^2 y)  \to   -\inf_{c > 0} c(\kappa(\mu) + J(y/c^2)) + \Lambda^*(x-c)
%\end{align}
%to establish \eqref{eq:uld_sparse_particulier}.

Notice that $\Lambda_0=\Lambda_X$. 
By \eqref{eq:maj_conditionnement} and Proposition \ref{prop:LLT_gd_dev},  making the change of variable $c=(y/\delta)^{1/2}$ in the infimum, and setting $x=(1-\mu)^{-1}$,  
it suffices to prove that 
\begin{align}\label{eq:ud_sparse_particulier}
\frac{1}{N} \log &\ \Prob(\rva{S}{N} = m_N,  \; \rva{T}{N} \geqslant N^2 y) \nonumber\\
 &\to   
\begin{cases}
 - \inf\limits_{c > 0} \big[c(\kappa(\mu) + J(y/c^2)) + \Lambda_X^*(x-c)\big]
& \text{if $y<(x-1)^2/2$}\\
-\infty & \text{if $y\geqslant (x-1)^2/2$}
\end{cases}
\end{align}
to establish \eqref{eq:uld_sparse_particulier}.
%\alert{
%\begin{remark} \leavevmode
%\begin{itemize}
%\item If $x-\sqrt{2y} > (1-\mu)^{-1}$, the formula may be different (see the end of the proof).
%\item $I(x,y) < \infty$ if and only if $y < (x-1)^2/2$. 
%\item What about the very sparse case $x = 1$?
%\end{itemize}
%\end{remark}
%}
Assume that \eqref{eq:ud_sparse_particulier} holds for $y < (x-1)^2/2$.
Observe that the probability in the left-hand side of \eqref{eq:ud_sparse_particulier} is decreasing in $y$. Moreover, if $y<(x-1)^2/2)$, then
\begin{align*}
\inf_{c > 0} \big[c(\kappa(\mu) + J(y/c^2)) + \Lambda_X^*(x-c)\big]
&= \inf_{\sqrt{2y}<c< x-1} \big[c(\kappa(\mu) + J(y/c^2)) + \Lambda_X^*(x-c)\big]\\
&\geqslant \inf_{\sqrt{2y}<c< x-1} c J(y/c^2)\\
&\geqslant \sqrt{2y} J(y/(x-1)^2)\\
&\to \infty \quad \text{as $y\to (x-1)^2/2$.}
\end{align*}
So \eqref{eq:ud_sparse_particulier} holds for $y\geqslant (x-1)^2/2$.

\begin{proof}[Proof of \eqref{eq:ud_sparse_particulier} -- Lower bound for $y< (x-1)^2/2$]
Let $c \in \intervalleoo{\sqrt{2y}}{x-1}$ and $l_N = \floor{cN}$. We have
\begin{align*}
\frac{1}{N} \log\ & \Prob(\rva{S}{N} = m_N,\ \rva{T}{N} \geqslant N^2 y) \\
 & \geqslant \frac{1}{N} \log \Prob(\rva{X}{N} = l_N, \rva{Y}{N} \geqslant N^2 y) + \frac{1}{N} \log \Prob(\rva[N-1]{S}{N} = m_N-l_N) \\
 &  =  - c(\kappa(\mu) + J(y/c^2)) - \Lambda_{\rva{X}{N}}^*((m_N - l_N) / N) +o(1),
\end{align*}
by Propositions \ref{prop:queue_couple} and \ref{prop:LLT_gd_dev}.
Now, let us prove that 
\begin{align}\label{eq:lambda0*}
\Lambda_{\rva{X}{N}}^*\left(\frac{m_N - l_N}{ N} \right) \to  \Lambda_X^*(x-c).
\end{align}
One has
\[
\Lambda_{\rva{X}{N}}^*\left(\frac{m_N - l_N}{ N} \right) = \sup_{s \in \R} \Bigl( s \frac{m_N - l_N}{N} - \Lambda_{\rva{X}{N}}(s) \Bigr) .
\]
The sequence of concave functions
\[
f_N(s) = s \frac{m_N - l_N}{N} - \Lambda_{\rva{X}{N}}(s)
\]
converges pointwise to the strictly concave function $f(s) = s(x - c) - \Lambda_X(s)$. Let $\tau$ be the unique point such that $f'(\tau) = 0$ (i.e\ $\Lambda_X'(\tau) = x - c$). Let $\varepsilon > 0$. Since $f'$ is decreasing, $f'(\tau - \varepsilon) > 0 > f'(\tau + \varepsilon)$. Now, the functions $f_N'$ converge to $f'$ uniformly on $\intervalleff{\tau - \varepsilon}{\tau + \varepsilon}$. So, for $N$ large enough, $f_N'(\tau - \varepsilon) > 0 > f_N'(\tau + \varepsilon)$. Therefore, for $N$ large enough, the supremum of $f_N$ is attained on $\intervalleff{\tau - \varepsilon}{\tau + \varepsilon}$ and converges to the supremum of $f$ (with $f_N(\tau_N) = \sup f_N$, $\sup f_N = f_N(\tau_N) \leqslant f(\tau_N) + \eta \leqslant \sup f + \eta$ and $\sup f_N \geqslant f_N(\tau) \geqslant f(\tau) - \eta$).
Hence,
\begin{align*}
\frac{1}{N} \log \Prob(\rva{S}{N} = m_N,\ \rva{T}{N} \geqslant N^2 y)
  \to  - c(\kappa(\mu) + J(y/c^2)) - \Lambda_X^*(x-c).
\end{align*}

Taking the supremum in $c > 0$ yields the desired lower bound.
\end{proof}

\begin{proof}[Proof of \eqref{eq:ud_sparse_particulier} -- Upper bound  for $y< (x-1)^2/2$]
Let us write
\begin{equation} \label{eq_PN0_PN1}
\Prob(\rva{S}{N} = m_N,\ \rva{T}{N} \geqslant N^2 y) = P_{N,0} + P_{N,1}
\end{equation}
where
\[
P_{N,0} = \Prob(\rva{S}{N} = m_N,\ \rva{T}{N} \geqslant N^2 y,\ \forall i \in \intervallentff{1}{N} \quad \rva[i]{Y}{N} \leqslant N^2 y)
\]
and
\[
P_{N,1} = \Prob(\rva{S}{N} = m_N,\ \rva{T}{N} \geqslant N^2 y,\ \exists i \in \intervallentff{1}{N} \quad \rva[i]{Y}{N} > N^2 y) .
\]

\paragraph{Behavior of $P_{N,0}$} 

Let us apply the exponential version of Chebyshev's inequality. Let $(s,t) \in \intervalleof{-\infty}{0} \times \intervallefo{0}{+\infty}$. We have
\begin{align}
P_{N,0}
 & \leqslant \Espe[\indic_{\rva{S}{N}-m_N \leqslant 0}\ \indic_{\rva{T}{N}/N-Ny \geqslant 0}\ \indic_{\forall i \in \intervallentff{1}{N} \quad \rva[i]{Y}{N} \leqslant N^2 y}] \nonumber \\
 & \leqslant \Espe[e^{s(\rva{S}{N}-m_N) + t(\rva{T}{N}/N-Ny)}\ \indic_{\forall i \in \intervallentff{1}{N} \quad \rva[i]{Y}{N} \leqslant N^2 y}] \nonumber \\
 & = e^{-N(s m_N/N + ty)} \Espe[e^{s\rva{X}{N} + t\rva{Y}{N}/N}\ \indic_{\rva{Y}{N} \leqslant N^2 y}]^N . \label{eq:PN0_maj}
\end{align}
Let us write
\begin{align}
\Espe[e^{s\rva{X}{N} + t\rva{Y}{N}/N}\ \indic_{\rva{Y}{N} \leqslant N^2 y}]
 & = \Espe[e^{s\rva{X}{N} + t\rva{Y}{N}/N}\ \indic_{\rva{Y}{N} \leqslant N^{1/2}}] + \Espe[e^{s\rva{X}{N} + t\rva{Y}{N}/N}\ \indic_{N^{1/2} < \rva{Y}{N} \leqslant N^2 y}] \nonumber \\
 & \eqdef E_1 + E_2 \label{eq:E1_E2_def} .
\end{align}

First, by uniform integrability (see Proposition \ref{prop:cv_moments}),
\begin{equation} \label{eq:E1_lim}
E_1  \to  \Espe[e^{sX}] .
\end{equation}

Secondly, remembering that $\Prob(\rva{X}{N}=l,\ \rva{Y}{N}=p) = 0$ if $l(l-1)/2 < p$,
\begin{align}
E_2
 & = \Espe[e^{s\rva{X}{N} + t\rva{Y}{N}/N}\ \indic_{N^{1/2} < \rva{Y}{N} \leqslant N^2 y}] \nonumber \\
 & = \sum_{\substack{l > N^{1/4} \\ N^{1/2} < p \leqslant N^2 y}} e^{sl+tp/N} \Prob(\rva{X}{N}=l,\ \rva{Y}{N}=p) \nonumber \\
 & = \sum_{l > N^{1/4}} e^{sl} \biggl( \sum_{N^{1/2} < p \leqslant N^2 y} e^{tp/N}(1-e^{-t/N}) \Prob(\rva{X}{N}=l,\ \rva{Y}{N} \geqslant p) \nonumber \\
 & \hspace{2cm} + \Prob(\rva{X}{N}=l,\ \rva{Y}{N} > N^{1/2}) e^{t\lfloor N^{1/2} \rfloor/N} - \Prob(\rva{X}{N}=l,\ \rva{Y}{N} > N^2 y) e^{t\lfloor N^2 y \rfloor / N} \biggr) , \label{eq:E2}
\end{align}
after a summation by parts. By uniform integrability (see Proposition \ref{prop:cv_moments}),
\begin{equation} \label{eq:E2_reste1}
\sum_{l > N^{1/4}} e^{sl} \Prob(\rva{X}{N}=l,\ \rva{Y}{N} > N^{1/2}) e^{t\lfloor N^{1/2} \rfloor/N} = \Espe[e^{s\rva{X}{N}+t\lfloor N^{1/2} \rfloor/N} \indic_{\rva{Y}{N} > N^{1/2}}]  \to  0 .
\end{equation}

The proof of the following lemma is postponed to the end of the paper (see page \pageref{proof_J_superconvexe}).

\begin{lem} \label{lem:J_superconvexe}
The function $K \colon \delta \in \intervallefo{0}{1/2} \mapsto \delta^{-1/2} J(\delta)$ is increasing, convex, and $K(\delta) \to \infty$ as $\delta \to 1/2$.
\end{lem}

For all $\varepsilon \in \intervalleoo{0}{1/2}$, we introduce the function
\begin{align*}
J_\varepsilon \colon \intervallefo{0}{\infty} & \to \R \\
\delta & \mapsto
\begin{cases}
J(\delta) & \text{if $\delta \leqslant 1/2-\varepsilon$} \\
\delta^{1/2} [ K(1/2-\varepsilon) + (\delta-1/2+\varepsilon) K'(1/2-\varepsilon) ] & \text{if $1/2-\varepsilon < \delta \leqslant 1/2 $} \\
\infty & \text{if $\delta > 1/2$.}
\end{cases}
\end{align*}

The following lemma is a straightforward consequence of Lemma \ref{lem:J_superconvexe}.

\begin{lem} \label{lem:Jepsilon}
The function $J_\varepsilon$ is non decreasing and less than $J$. Moreover, the function $\delta \in \intervalleff{0}{1/2} \mapsto \delta^{-1/2} J_\varepsilon(\delta)$ is a bounded convex function.
\end{lem}

Let $\varepsilon \in \intervalleoo{0}{\kappa(\mu) \wedge 1/2}$. As a consequence of \eqref{eq:E2}, \eqref{eq:E2_reste1}, and Lemmas \ref{lem:maj_queue_XY_uniforme} and \ref{lem:Jepsilon}, we get, for any $a > 0$,
\begin{align}
\limsup E_2
 & \leqslant \limsup \sum_{l > N^{1/4}} e^{sl} \sum_{N^{1/2} < p \leqslant N^2 y} e^{tp/N}(1-e^{-t/N}) e^{-l(\kappa(\mu) + J_\varepsilon(p/l^2) - \varepsilon)} \nonumber \\
  & \leqslant \limsup \Biggl( \sum_{\substack{N^{1/4} < l \leqslant Na \\ N^{1/2} < p \leqslant N^2 y}} e^{l(s-\kappa(\mu) - J_\varepsilon(p/l^2) + \varepsilon)+tp/N} \nonumber + N^2 y \sum_{l> Na} e^{-(\kappa(\mu)-\varepsilon)l+tNy} \Biggr) \nonumber \\
  & \leqslant \limsup N^3 ay \cdot \exp\Biggl( N \cdot \max_{\substack{N^{1/4} <l \leqslant Na \\ N^{1/2} < p \leqslant N^2 y}} \Bigl[ \frac{l}{N}\Bigl(s-\kappa(\mu) - J_\varepsilon\Bigl(\frac{p}{l^2}\Bigr)  + \varepsilon \Bigr)+\frac{tp}{N^2} \Bigr] \Biggr), \label{eq:E2_max}
\end{align}
as soon as $(\kappa(\mu)-\varepsilon)a > ty$. Now,
\begin{align}
\max_{\substack{N^{1/4} < l \leqslant Na \\ N^{1/2} < p \leqslant N^2 y}} & \Bigl[ \frac{l}{N}\Bigl(s-\kappa(\mu) - J_\varepsilon\Bigl(\frac{p}{l^2}\Bigr) + \varepsilon \Bigr)+\frac{tp}{N^2} \Bigr] \nonumber\\
 & \leqslant \sup_{\substack{0 < \delta < 1/2 \\ N^{-3/4} < u \leqslant \sqrt{y\delta^{-1}}}} \bigl[ u(s-\kappa(\mu) - J_\varepsilon(\delta) + \varepsilon)+t\delta u^2 \bigr] 
\eqdef S . \label{eq:S_def}
\end{align}
As soon as
\[
\sup_{0 < \delta < 1/2} \bigl[ \sqrt{y \delta^{-1}} (s-\kappa(\mu) - J_\varepsilon(\delta) + \varepsilon) + ty \bigr] < 0 ,
\]
i.e.\ as soon as
\begin{equation} \label{eq:hyp_t}
t < \inf_{0 < \delta < 1/2} \frac{\kappa(\mu)+J_\varepsilon(\delta)-\varepsilon-s}{\sqrt{\delta y}} ,
\end{equation}
and, for all $N$ large enough,
\begin{align}
S
 & = \sup_{0 < \delta < 1/2} \bigl[ N^{-3/4} (s - \kappa(\mu) - J_\varepsilon(\delta) + \varepsilon) + t\delta N^{-3/2} \bigr]
 \leqslant -\frac{(\kappa(\mu)-\varepsilon) N^{-3/4}}{2} . \label{eq:S_maj}
\end{align}
Therefore, under \eqref{eq:hyp_t}, using \eqref{eq:E2_max}, \eqref{eq:S_def}, and \eqref{eq:S_maj},
\begin{equation} \label{eq:E2_lim}
\limsup E_2 \leqslant \limsup N^3 ay \cdot \exp\Bigl( -\frac{(\kappa(\mu)-\varepsilon) N^{1/4}}{2} \Bigr) = 0 .
\end{equation}
Combining \eqref{eq:PN0_maj}, \eqref{eq:E1_E2_def}, \eqref{eq:E1_lim}, and \eqref{eq:E2_lim}, we get
\begin{align*}
\limsup \frac{1}{N} \log P_{N,0}
 & \leqslant \inf\left\{ \Lambda_X(s) - sx - ty\ ;\ s \leqslant 0,\ t < \inf_{0 < \delta < 1/2} (\kappa(\mu)+J_\varepsilon(\delta)-\varepsilon-s)(\delta y)^{-1/2} \right\} \\
 & = \inf_{s \leqslant 0} \sup_{0 < \delta < 1/2} \Lambda_X(s) - sx - (\kappa(\mu)+J_\varepsilon(\delta)-\varepsilon-s)(y/\delta )^{1/2}  \\
 & \eqdef M_\varepsilon' .
\end{align*}

The proof of the following lemma is postponed to the end of the paper (see page \pageref{proof_concave_convexe}).

\begin{lem}\label{lem:concave_convexe}
The function
\begin{align*}
f \colon \intervalleoo{0}{1/2} \times \intervalleof{-\infty}{0} & \to \R \\
(\delta, s) & \mapsto \Lambda_X(s) - sx - (\kappa(\mu)+J_\varepsilon(\delta)-\varepsilon-s)(y/\delta )^{1/2} 
\end{align*}
is concave in $\delta$ and convex in $s$. Moreover, $f(\delta_0, \cdot)$ is bounded from below for some $\delta_0 \in \intervalleoo{0}{1/2}$.
\end{lem}

Thanks to Lemma \ref{lem:concave_convexe}, the minimax theorem of \cite{2015_PerchetVigneral_AMinmaxTheorem} applies and yields
\begin{align*}
M_\varepsilon'
 & = \sup_{0 < \delta < 1/2} \inf_{s \leqslant 0} \Bigl[ \Lambda_X(s) - sx - (\kappa(\mu)+J_\varepsilon(\delta)-\varepsilon-s) (y/\delta )^{1/2} \Bigr] \\
 & = - \inf_{0 < \delta < 1/2} \Bigl[ (y/\delta )^{1/2} (\kappa(\mu)+J_\varepsilon(\delta)-\varepsilon) + \Lambda_X^*\bigl( x - (y/\delta )^{1/2} \bigr) \Bigr] .
\end{align*}
Notice that, by Lemma \ref{lem:J_superconvexe} and since $x - (y/\delta )^{1/2} \leqslant x - \sqrt{2y} \leqslant (1-\mu)^{-1}$,
\begin{align*}
\inf_{1/2-\varepsilon < \delta < 1/2} &\Bigl[ (y/\delta )^{1/2} (\kappa(\mu)+J_\varepsilon(\delta)-\varepsilon) + \Lambda_X^*\bigl( x - (y/\delta )^{1/2} \bigr) \Bigr] \\
&\geqslant \sqrt{2y} (\kappa(\mu)-\varepsilon) + \sqrt{y}K(1/2-\varepsilon)+\Lambda_X^*\bigl( x-\sqrt{2y} \bigr)\xrightarrow[\varepsilon \to 0]{} \infty .
\end{align*}
A fortiori, since $J_\varepsilon \leqslant J$,
\[
\inf_{1/2-\varepsilon < \delta < 1/2} \Bigl[ (y/\delta )^{1/2} (\kappa(\mu)+J(\delta)-\varepsilon) + \Lambda_X^*\bigl( x - (y/\delta )^{1/2} \bigr) \Bigr] \xrightarrow[\varepsilon \to 0]{} \infty .
\]
So, if $\varepsilon$ is small enough, i.e.\ $\varepsilon \in \intervalleoo{0}{\varepsilon_0}$,
\[
M_\varepsilon' = - \inf_{0 < \delta < 1/2} \Bigl[ (y/\delta )^{1/2} (\kappa(\mu)+J(\delta)-\varepsilon) + \Lambda_X^*\bigl( x - (y/\delta )^{1/2} \bigr) \Bigr] .
\]
Finally, again applying the minimax theorem of \cite{2015_PerchetVigneral_AMinmaxTheorem}, we get
\begin{align*}
\inf_{0 < \varepsilon < \varepsilon_0} M_\varepsilon'
 & = - \sup_{0 < \varepsilon < \varepsilon_0} \inf_{0 < \delta < 1/2} \Bigl[ (y/\delta )^{1/2} (\kappa(\mu)+J(\delta)-\varepsilon) + \Lambda_X^*\bigl( x - (y/\delta )^{1/2} \bigr) \Bigr] \\
 & = - \sup_{0 < \varepsilon < \varepsilon_0} \inf_{\sqrt{2y} < c \leqslant x-1} \Bigl[ c (\kappa(\mu)+J(y / c^2)-\varepsilon) + \Lambda_X^*(x - c) \Bigr] \\
% & = - \inf_{c > \sqrt{2y}} \Bigl[ c (\kappa(\mu)+J(y / c^2)-\varepsilon) + \Lambda_X^*(x - c) \Bigr] \\
 & = - \inf_{\sqrt{2y} < c \leqslant x-1} \Bigl[ c(\kappa(\mu)+J(y/c^2)) + \Lambda_X^*(x-c) \Bigr] ,
\end{align*}
since $\Lambda_X^*(x-c) = \infty$ for $c > x-1$ and since the function
\begin{align*}
g \colon \intervalleoo{0}{\varepsilon_0} \times \intervalleof{\sqrt{2y}}{x-1} & \to \R \\
(\varepsilon, c) & \mapsto c (\kappa(\mu)+J(y / c^2)-\varepsilon) + \Lambda_X^*(x - c)
\end{align*}
is nonnegative, concave in $\varepsilon$ and convex in $c$ (note that $c \mapsto c J(y / c^2)$ is convex by differentiating twice and applying Lemma \ref{lem:J_superconvexe}). So, we have proved that
\begin{equation} \label{eq_PN0_maj}
\limsup \frac{1}{N} \log P_{N,0} \leqslant - \inf\limits_{c > 0} \big[c(\kappa(\mu) + J(y/c^2)) + \Lambda_X^*(x-c)\big] .
\end{equation}

\medskip

\paragraph{Behavior of $P_{N,1}$}

Let $\varepsilon > 0$. We have
\begin{align*}
& P_{N,1}
  = \Prob(\rva{S}{N} = m_N,\ \rva{T}{N} \geqslant N^2 y,\ \exists i \in \intervallentff{1}{N} \quad \rva[i]{Y}{N} > N^2 y) \\
 & \leqslant N \Prob(\rva{S}{N} = m_N,\ \rva[N]{Y}{N} > N^2 y) \\
 & \leqslant N \sum_{l=1}^\infty \Prob\Bigl( m_N-Nl\varepsilon \leqslant \rva[N-1]{S}{N} < m_N-N(l-1)\varepsilon \Bigr) \Prob( N(l-1)\varepsilon < \rva{X}{N} \leqslant Nl\varepsilon,\ \rva{Y}{N} > N^2 y) \\
 & \leqslant N \sum_{l=1}^\infty P_{N,1,l} ,
\end{align*}
where
\[
P_{N,1,l} \defeq \Prob(\rva[N-1]{S}{N} < m_N-N(l-1)\varepsilon) \Prob( N(l-1)\varepsilon < \rva{X}{N} \leqslant Nl\varepsilon,\ \rva{Y}{N} > N^2 y) .
\]
Now, if $m_N/N - (l-1) \varepsilon \leqslant 0$, then $P_{N, 1, l} = 0$. Let $L \defeq \sup \ensavec{(1 + (m_N / N) / \varepsilon)}{N \geqslant 1} < \infty$ since $m_N / N \to x$. Consequently,
\begin{align*}
P_{N,1} & \leqslant N \sum_{l=1}^L P_{N,1,l} .
\end{align*}
%$K$ is the smallest integer such that $(K\varepsilon)^2/2 \geqslant y$ (remember that $\rva{Y}{N} \leqslant \rva{X}{N}^2/2$), and $L$ is a fixed integer such that, for all $N \geqslant 1$, $m_N/N - (L-1)\varepsilon \leqslant 0$ (remember that $m_N/N \to x$).
Let us evaluate each term. For all $l \in \intervallentff{1}{L}$, 
\begin{align*}
\limsup \frac{1}{N} \log P_{N,1,l}
 & \leqslant -\Lambda_X^*(x-(l-1)\varepsilon) - l\varepsilon \Bigl(\kappa(\mu)+J\Bigl(\frac{y}{(l\varepsilon)^2}\Bigr)\Bigr) .
\end{align*}
applying the exponential version of Chebyshev's inequality, \eqref{eq:lambda0*}, and Proposition \ref{prop:queue_couple}. 
Applying the principle of the largest term, we get:
\[
\limsup \frac{1}{N} \log P_{N,1} \leqslant - \inf_{c > \sqrt{2y}} \left[ \Lambda_X^*(x-c+\varepsilon) + c \Bigl(\kappa(\mu)+J\Bigl(\frac{y}{c^2}\Bigr)\Bigr) \right],
\]
and we get the desired upper bound for $P_{N,1}$ when $\varepsilon \to 0$, since $\Lambda_X^*(a) = \infty$ for $a < 1$, since $\Lambda_X^*(x-c+\varepsilon)$ converges uniformly in $c \in \intervalleof{\sqrt{2y}}{x-1}$ to $\Lambda_X^*(x-c)$, and since
\begin{align*}
\inf_{x-1 < c \leqslant x-1+\varepsilon} \left[\Lambda_X^*(x-c+\varepsilon) + c \Bigl(\kappa(\mu)+J\Bigl(\frac{y}{c^2}\Bigr)\Bigr)\right]
 & \geqslant \Lambda_X^*(1+\varepsilon) + \inf_{x-1 < c \leqslant x-1+\varepsilon} c \Bigl(\kappa(\mu)+J\Bigl(\frac{y}{c^2}\Bigr)\Bigr) \\
 & \xrightarrow[\varepsilon \to 0]{} \Lambda_X^*(1) + (x-1) \Bigl(\kappa(\mu)+J\Bigl(\frac{y}{(x-1)^2}\Bigr)\Bigr). 
\end{align*}
So, we have proved that
\begin{equation} \label{eq_PN1_maj}
\limsup \frac{1}{N} \log P_{N,1} \leqslant - \inf\limits_{c > 0} \big[c(\kappa(\mu) + J(y/c^2)) + \Lambda_X^*(x-c)\big] .
\end{equation}
\end{proof}

The proof of \eqref{eq:ud_sparse_particulier} follows from \eqref{eq_PN0_PN1}, \eqref{eq_PN0_maj}, and \eqref{eq_PN1_maj} and the principle of the largest term.
%Let us turn to the proof of the remaining technical Lemmas.
%Remember that, for all $\delta \in \intervallefo{0}{1/2}$,
%\[
%J(\delta) = \Bigl( \frac{1}{2} - \delta \Bigr)  \lambda(\delta) + \log\Bigl( 1 -  \Bigl( \frac{1}{2} + \delta \Bigr)  \lambda(\delta) \Bigr)
%\]
%where $\lambda(\delta)$ is the smallest solution of the equation in $\lambda$
%\[
%\Bigl(   \Bigl( \frac{1}{2} + \delta \Bigr) \lambda(\delta) - 1 \Bigr)(1-e^{\lambda(\delta)}) = \lambda(\delta) .
%\]
It remains to prove Lemmas \ref{lem:J_superconvexe} and \ref{lem:concave_convexe}.  Let us define the function
\begin{align*}
\delta \colon \intervalleof{-\infty}{0} & \to \intervallefo{0}{1/2} \\
\lambda & \mapsto
\begin{cases}
0 & \text{if $\lambda = 0$} \\
\frac{1}{\lambda} + \frac{1}{1-e^\lambda} - \frac{1}{2} & \text{if $\lambda \in \intervalleoo{-\infty}{0}$.}
\end{cases}
\end{align*}

An easy computation shows that the function $\delta$ is a smooth nonnegative and concave decreasing bijection. The function $\lambda \colon \intervallefo{0}{1/2} \to \intervalleof{-\infty}{0}$ defined by \eqref{eq:lambda} is the inverse bijection of $\delta$ and thus is a smooth decreasing function. 
Now let us introduce the functions
\begin{align*}
H \colon \intervalleof{-\infty}{0} & \to \intervallefo{0}{+\infty} \\
\lambda & \mapsto \lambda \cdot \Bigl( \frac{1}{2} - \delta(\lambda) \Bigr) + \log\Bigl( 1 - \lambda \cdot \Bigl( \frac{1}{2} + \delta(\lambda) \Bigr) \Bigr)
\end{align*}
and
\begin{align*}
F \colon \intervalleof{-\infty}{0} & \to \intervallefo{0}{+\infty} \\
\lambda & \mapsto
\begin{cases}
0 & \text{if $\lambda = 0$} \\
\delta(\lambda)^{-1/2} H(\lambda) & \text{if $\lambda \in \intervalleoo{-\infty}{0}$.}
\end{cases}
\end{align*}

\begin{proof}[Proof of Lemma \ref{lem:J_superconvexe}] \label{proof_J_superconvexe}
The fact that $K(\delta) \to \infty$ as $\delta \to 1/2$ follows from the already mentioned fact that $J(\delta) \to \infty$ as $\delta \to 1/2$. We want to prove that $\delta \mapsto K(\delta) = \delta^{-1/2} J(\delta) = F(\lambda(\delta))$ is an increasing convex function. Since $K(\delta) \sim 6\delta^{3/2}$ as $\delta \to 0$, $K'(0) = 0$, so it only remains to prove that $K'' > 0$ over $\intervalleoo{0}{1/2}$. Using the expressions of the first and second derivatives of the inverse function $\lambda$ of $\delta$, one gets:
\begin{align*}
K''(\delta)
 & = \frac{d^2}{d\delta^2} F(\lambda(\delta)) \\
 & =\lambda''(\delta)F'(\lambda(\delta))+(\lambda'(\delta))^2 F''(\lambda(\delta))\\
 & = \frac{1}{\delta'(\lambda(\delta))}\left(-\frac{\delta''(\lambda(\delta))}{\delta'(\lambda(\delta))^2}F'(\lambda(\delta))+\frac{1}{\delta'(\lambda(\delta))}F''(\lambda(\delta))\right)\\
 & =  \frac{1}{\delta'(\lambda(\delta))}\left(\frac{F'}{\delta'}\right)'(\lambda(\delta)).
\end{align*}
Hence, since $\delta' < 0$, our study reduces to show that $F'/\delta'$ is a decreasing function over $\intervalleoo{-\infty}{0}$. Now, straightforward calculations yield the magical identity $H'(\lambda)=-\lambda \delta'(\lambda)$. Hence
\begin{align*}
\frac{F'(\lambda)}{\delta'(\lambda)}
 & = - \delta(\lambda)^{1/2} \cdot \frac{\lambda \delta(\lambda)+\frac 12 H(\lambda)}{\delta(\lambda)^2} \eqdef - \delta(\lambda)^{1/2} f(\lambda) .
\end{align*}
and also
\begin{align*}
f'(\lambda) = \frac{\delta(\lambda)^2-\frac 32 \lambda \delta'(\lambda) \delta(\lambda)-H(\lambda)\delta'(\lambda)}{\delta(\lambda)^3} \eqdef \frac{k(\lambda)}{\delta(\lambda)^3}.
\end{align*}
Differentiating the function $k$, we get
\begin{align*}
k'(\lambda)=\frac 12 \delta'(\lambda)\bigl(\delta(\lambda)-\lambda \delta'(\lambda)\bigr) - \delta''(\lambda)\Bigl(H(\lambda)+\frac 32 \lambda \delta(\lambda)\Bigr) .
\end{align*}

On the one hand, $H(\lambda) + \frac 32 \lambda \delta(\lambda) < 0$, because
\[
\frac{d}{d\lambda} \Bigl( H(\lambda) + \frac 32 \lambda \delta(\lambda) \Bigr) = \frac{e^{2\lambda}(-3\lambda+4)+e^\lambda(2\lambda^2-8)+3\lambda+4}{4\lambda(1-e^\lambda)^2} ,
\]
the sign of which is easy to find (by differentiating several times). On the other hand, $\delta(\lambda)-\lambda \delta'(\lambda) > 0$, because
\[
\delta(\lambda)-\lambda \delta'(\lambda)
 % = \frac 2\lambda -\frac 12 +\frac{1-(\lambda+1)e^\lambda }{(1-e^\lambda)^2}
 % = \frac{2\lambda (1-(\lambda+1)e^\lambda) +(4-\lambda)(1-e^\lambda)^2}{2\lambda (1-e^\lambda)^2}
 = \frac{e^{2\lambda}(-\lambda+4) + e^\lambda(-2\lambda^2-8) + \lambda+4}{2\lambda (1-e^\lambda)^2} ,
\]
the sign of which is easy to find (similarly). Henceforth, $k' < 0$ and $k > 0$ over $\intervalleoo{-\infty}{0}$ ($k$ is decreasing on $\intervalleof{-\infty}{0}$ and $k(0)=0$). Finally, $f$ is an increasing and nonpositive function. Together with the fact that $\lambda \mapsto -\delta(\lambda)^{1/2}$ is also an increasing and nonpositive function, we finally get that $F'/\delta'$ is a decreasing function.
\end{proof}

\begin{proof}[Proof of Lemma \ref{lem:concave_convexe}] \label{proof_concave_convexe}
The concavity of $f(\cdot, s)$ follows from Lemma \ref{lem:Jepsilon} and the fact that $\kappa(\mu) - \varepsilon - s \geqslant 0$ and $y \geqslant 0$. The convexity of $f(\delta,\cdot)$ follows from the convexity of $\Lambda_X$. Finally, since $y < (x-1)^2/2$, one can choose $\delta_0 \in \intervalleoo{0}{1/2}$ such that $x-(y/\delta_0)^{1/2} > 1$. Finally,
\[
\Lambda_X(s)
 = \log \Espe[e^{sX}]
 \geqslant \log(e^s \Prob(X = 1))
 = s + \log\Prob(X=1) ,
\]
so the function $f(\delta_0,\cdot)$ is bounded from below.
\end{proof}

%\begin{proof}[Proof of Lemma \ref{lem:delta}]
%The function $\delta$ is continuous over $\intervalleof{-\infty}{0}$. Over $\intervalleoo{-\infty}{0}$, $\delta$ is twice differentiable and one has
%\[
%\delta'(\lambda) = \frac{-e^{2\lambda} + (\lambda^2+2)e^\lambda - 1}{\lambda^2(1-e^\lambda)^2}
%\]
%and
%\[
%\delta''(\lambda)
% % = \frac{\lambda^3 e^\lambda(1+e^\lambda)+ 2(1-e^\lambda)^3}{\lambda^3(1-e^\lambda)^3} 
% = \frac{-2e^{3\lambda} + e^{2\lambda}(\lambda^3+6) + e^\lambda(\lambda^3-6) + 2}{\lambda^3(1-e^\lambda)^3}
%.
%\]
%The latter denominator is negative and a study of the numerator shows that it is positive. Hence, $\delta$ is concave. Moreover, $\delta'$ is decreasing and  
%since $\delta'(\lambda) \sim -\lambda^{-2} < 0$ as $\lambda \to -\infty$, $\delta'$ is negative. Therefore, $\delta$ is decreasing. The remaining conclusions of the lemma are obvious. 
%%
%%\begin{center}
%%\begin{tikzpicture}
%%   \tkzTabInit{$\lambda$ / 1 , $\delta''(\lambda)$ / 1, $\delta'(\lambda)$ / 1, $\delta(\lambda)$ / 1.5}{$-\infty$, $0$}
%%   \tkzTabLine{, -, z }
%%   \tkzTabLine{, -,  }
%%   \tkzTabVar{+/ ${1/2}$, -/ $0$}
%%\end{tikzpicture}
%%\end{center}
%\end{proof}

\paragraph{Acknowledgment} We deeply thank the anonymous reviewer for his thorough reading of our manuscript and for his insightful comments. 

\bibliographystyle{abbrv}
\bibliography{biblio_gde_dev_alpha}

\end{document}